\documentclass[12pt]{amsart}
\usepackage{amssymb}
\usepackage{amscd}
\usepackage{amsmath}

\newcommand{\ga}{{\mathfrak a}}
\newcommand{\gb}{{\mathfrak b}}
\newcommand{\gog}{{\mathfrak g}}
\newcommand{\gh}{{\mathfrak h}}
\newcommand{\gl}{{\mathfrak l}}
\newcommand{\gm}{{\mathfrak m}}
\newcommand{\gn}{{\mathfrak n}}
\newcommand{\gp}{{\mathfrak p}}
\newcommand{\gs}{{\mathfrak s}}
\newcommand{\bB}{{\bf B}}
\newcommand{\bG}{{\bf G}}

\newcommand{\bL}{{\bf L}}
\newcommand{\bM}{{\bf M}}

\newcommand{\bP}{{\bf P}}
\newcommand{\bS}{{\bf S}}
\newcommand{\bT}{{\bf T}}

\newcommand{\bD}{{\bf D}}

\newcommand{\al}{\alpha}
\newcommand{\be}{\beta}

\newcommand{\eps}{\epsilon}

\newcommand{\cal}{\mathcal} 
\newcommand{\Cscr}{{\cal C}}
\newcommand{\Nscr}{{\cal N}}

\newcommand{\Oscr}{{\cal O}}
\newcommand{\Vscr}{{\cal V}}

\newcommand{\Tscr}{{\cal T}}
\newcommand{\Iscr}{{\cal I}}
\newcommand{\Jscr}{{\cal J}}
\newcommand{\Wscr}{{\cal W}}

\newcommand{\Mscr}{{\cal M}}

\newcommand{\Xscr}{{\cal X}}
\newcommand{\Lie}{{\rm Lie\,}}
\newcommand{\sh}{{\rm sh\,}}

\newcommand{\gr}{{\rm gr\,}}

\newcommand{\Ann}{{\rm Ann\,}}
\newcommand{\Spa}{{\rm span\,}}
\newcommand{\QED}{\par \hspace{15cm}$\blacksquare$ \par}
\newcommand{\pr}{^{\prime}}

\newcommand{\st}{\subset}
\newcommand{\Co}{{\mathbb C}} 

\newcommand{\Pf}{\noindent{\bf Proof.}\par\noindent}

\newcommand{\sr}{\scriptscriptstyle}
\newcommand{\vb}{\vrule height 14pt depth 7pt}
\newcommand{\ts}{\tabskip 4pt}
\newcommand{\vsa}{\noalign{\vskip-7pt}}
\newcommand{\ssa}{\noalign{\vskip -1pt}}

\newcommand{\bk}{\break}
\newcommand{\ov}{\overline}
\newcommand{\rar}{\rightarrow}

\newcommand{\Da}{\Downarrow}
\newcommand{\Ra}{\Rightarrow}

\newcommand{\ha}{\hookrightarrow}

\newcommand{\parno}{\par\noindent}
\newcommand{\dor}{\stackrel{\rm D}{\leq}}     
\newcommand{\dg}{\stackrel{\rm D}{\geq}}
\newcommand{\dos}{\stackrel {\rm D}{<}}
\newcommand{\dgs}{\stackrel {\rm D}{>}}
\newcommand{\dvo}{\stackrel{\rm D-V}{\leq}}     

\newcommand{\dvonot}{\stackrel{\rm D-V}{\not\leq}}
\newcommand{\go}{\stackrel {\rm G}{\leq}}     

\newcommand{\gegs}{\stackrel {\rm G}{>}}
\newcommand{\gos}{\stackrel {\rm G}{<}}
\newcommand{\cho}{\stackrel {\rm  Ch}{\leq}}     
\newcommand{\chos}{\stackrel {\rm  Ch}{<}}

\newcommand{\chg}{\stackrel {\rm  Ch}{\geq}}
\newcommand{\chvo}{\stackrel {\rm  V-Ch}{\leq}}     
\newcommand{\chvos}{\stackrel {\rm  V-Ch}{<}}

\newcommand{\ao}{\stackrel {\rm  A}{\leq}}     
\newcommand{\aos}{\stackrel {\rm  A}{<}}

\newcommand{\ag}{\stackrel {\rm  A}{\geq}}
\newtheorem*{theorem}{Theorem}
\newtheorem*{lemma}{Lemma}
\newtheorem*{defi}{Definition}
\newtheorem*{cor}{Corollary}
\newtheorem*{prop}{Proposition}
\newtheorem*{conj}{Conjecture}
\begin{document}
\title[orbital varieties]
{\bf On orbital variety closures in $\mathfrak{sl}_n$\\
III.Geometric properties}
\author{Anna Melnikov}
\address{Department of Mathematics,
University of Haifa,
Haifa 31905, Israel}
\email{melnikov@math.haifa.ac.il}
\begin{abstract}
This is the third paper in the series. Here we define a few
combinatorial orders on Young tableaux. The first order is obtained
from induced Duflo order by  the extension with the help of Vogan
$\Tscr_{\al, \be}$ procedure. We call it Duflo-Vogan order. The
second order is obtained from the generalization of Spaltenstein's
construction by consideration of an orbital variety as a double
chain of nilpotent orbits. We call it the chain order. Again,  we
use Vogan's $\Tscr_{\al, \be}$ procedure, however, this time to
restrict the chain order. We call it Vogan-chain order. The order on
Young tableaux defined by the inclusion of orbital variety closures
is called a geometric order and the order on Young tableaux defined
by inverse inclusion of primitive ideals is called an algebraic
order.

We get the following relations between the orders:
Duflo-Vogan order is an extension of the induced Duflo order;
the algebraic order is an extension of Duflo-Vogan order;
the geometric order is an extension of the algebraic order;
Vogan-chain order is an extension of the geometric order;
and, finally, the chain order is an extension of Vogan-chain order. The computations
show that Duflo-Vogan and Vogan-chain orders coincide on $\gs\gl_n$
for $n\leq 9$ and in $n=10$ there is one case (up to $\Tscr_{\al\be}$ procedure and
transposition) where Vogan-chain
order is a proper extension  of Duflo-Vogan order. In this only case
the algebraic order coincides with Vogan-chain order. These computations
permit us to conjecture that in $\gs\gl_n$ the algebraic order coincides
with the geometric order. As well we conjecture that the combinatorics of both
the inclusions on primitive ideals and  on orbital variety closures is defined
by Vogan-chain order on Young tableaux.
\end{abstract}
\maketitle

\section {Introduction}
\subsection{}\label{1.1}
This is the third paper in the series of three papers. We refer to the first two papers \cite{M4},
\cite{M5} as
Part I and Part II respectively. Our main objects in these series are orbital variety closures in
$\gs\gl_n.$ They are parameterized by Young tableaux. The purpose is to construct the combinatorial
order on Young tableaux defined in terms of Young tableaux only, describing inclusions of these closures.
We call this order a geometric order on Young tableaux.
\par
We begin with the description of the connection between orbital varieties in a semisimple Lie algebra
$\gog$ and primitive ideals in its enveloping algebra $U(\gog),$  containing the augmentation ideal
of the centre $Z(\gog)$ of $U(\gog).$ The role of orbital varieties in study of primitive
ideals was described in short in Part I, 1.3, however, here we would like to consider the connection
between these objects in more detail since
on one hand, the theory of  primitive ideals was the source of our interest to orbital varieties,
and on the other hand, the methods invented for the study
of primitive ideals can be successfully implemented to the study of orbital varieties,
especially, in the case of $\gog=\gs\gl_n.$

\subsection{}\label{1.2}
Let us set up the notation.
Let $\bG$ be a connected simply-connected complex algebraic group. Set $\gog=\Lie(\bG)$
and let $U(\gog)$ be the enveloping algebra of $\gog.$
Consider a co-adjoint action of $\bG$ on $\gog^*.$ Identify $\gog^*$ with $\gog$
through the Killing form. A $\bG$ orbit $\Oscr$ in $\gog$ is called nilpotent
if it consists of ad-nilpotent elements.
\par
Fix a triangular decomposition $\gog=\gn^-\oplus\gh\oplus\gn.$
Let $W$ be the Weyl group of $<\gog, \gh>.$
\par
Let $\Oscr$ be some nilpotent orbit. An irreducible component of
$\Oscr\cap\gn$
is called an orbital variety associated to $\Oscr.$ Recall from Part I, 2.1.3 that there exists
a surjection from $W$ onto the set of orbital varieties  defined by Steinberg's construction.
( Explicitly, let $\bB$ be the Borel subgroup of $\bG$ with $\Lie(\bB)=\gh\oplus \gn.$ Let
$\bB$ act adjointly on $\gn.$ Then each orbital variety closure $\ov\Vscr$ is
$\ov\Vscr_w=\ov{\bB(\gn\cap^w\gn)}$ for some $w\in W$). The fibres of this surjection
are called geometric cells.

\subsection{}\label{1.3}
Let $\Xscr_0$ denote the set of primitive ideals of
$U(\gog)$ containing the augmentation ideal of the
centre $Z(\gog)$ of $U(\gog).$ After M. Duflo ~\cite{D},
there exists a surjective map $\psi:W\rar \Xscr_0$ whose
fibres are called the algebraic (left) cells of $W$ (cf. \ref{4.1})
The inclusion relation on the primitive ideals
gives a partial order relation
on the left cells. We call it an algebraic order.
Its form in terms of the multiplicities in the
composition series of principal series
representations of $\bG$ was conjectured
by A. Joseph \cite{Jo}
and was shortly afterwards established by
D. Vogan \cite{V1}.
This result was later made purely combinatorial
by D. Kazhdan and G. Lusztig \cite{K-L}. Respectively, the algebraic order
is called also Kazhdan-Lusztig order in the literature.

\subsection{}\label{1.4}
As shown in \cite{B-B I}
and \cite{Jo1}, the associated variety of a primitive ideal
is a nilpotent orbit and, thus,
the Duflo map $\psi$ gives rise to a map from
$W$ to the set of nilpotent orbits.
However, this is generally not surjective.
The orbits in the image of this map are called Lusztig's special orbits.
Thus, despite the optimistic predictions of the orbit
method, it turns out that at our present level of
refinement
geometry of orbital varieties differs slightly
from representation theory of the corresponding Lie
algebras.
\par
The above considerations can be refined using the
associated variety
of a simple highest weight module.
As shown in \cite{B-B} and \cite{Jo1}, an irreducible component of such an associated variety
is the closure of some orbital variety.
Moreover, as
shown in \cite{B-B}, an inclusion
of primitive ideals implies the reverse inclusion
of corresponding associated varieties.

\subsection{}\label{1.5}
Let us explain the connection between
primitive ideals and orbital varieties in
terms of Goldi rank polynomials.

Let $R\st\gh^*$ denote the set of non-zero roots,  $R^+$ the set
of positive roots corresponding to $\gn$ in the triangular
decomposition of $\gog,$ and $\Pi\st R^+$ the resulting set of
simple roots. Set $\rho=0.5\sum\limits_{\alpha\in R^+}\alpha.$

For $w\in W$ let $L_w$ denote a simple highest weight
module with the highest weight $-w(\rho)-\rho.$
The formal character
of $L_w$ provides a polynomial $p_w$ on $\gh^*$
which by \cite{Jo2} determines the Goldie rank of the
corresponding primitive quotient and is called respectively
a Goldi rank polynomial.
A. Joseph further attached a characteristic
polynomial $q_{w}$ to an orbital variety
$\Vscr_w$ (cf. \cite{Jo1}). Now $p_{w^{\sr-1}}$
determines the characteristic polynomial of the
associated variety of the simple highest weight
module.
\par
The relation between geometric
cells defined in \ref{1.2} and algebraic  cells defined in \ref{1.3}
can be expressed in
terms of relation between $p_{w^{\sr-1}}$ and  $q_w.$
Unfortunately the difference between geometric
picture coming from Steinberg's construction
and the picture coming from primitive ideals
is, somehow, responsible for different complications
such as existence of special and
non-special orbits, mentioned in \ref{1.4}, and
reducibility (in general) of associated varieties of the
simple highest weight module. In particular,
the relationship between algebraic and
geometric cells is rather complicated, so that
it is not true that algebraic cell is a union
of corresponding geometric cells.

\subsection{}\label{1.6}
For $\gog=\gs\gl_n$ the above simplifies
considerably. Here all the orbits are special.
Moreover, as shown in \cite{M1}, the associated variety of a
simple highest weight module is always irreducible.
In particular, this result determines  the characteristic
polynomial of an orbital variety to be
$p_{w^{\sr-1}}$ for some $w\in W.$ Again up to interchanging
$w$ and $w^{\sr -1}$ geometric and algebraic cells coincide
and are further given by
Robinson-Schensted algorithm. This was first
observed by A.Joseph in the primitive ideal framework
and then by N. Spaltenstein and by R. Steinberg in the framework of orbital varieties.
\par
Let us denote by $\bT_n$ the set of standard Young tableaux with $n$
boxes. Thus, for  $T\in \bT_n$ we can uniquely define a primitive
ideal $I_T$ and an orbital variety $\Vscr_T.$ Correspondingly, we
define an algebraic order on Young tableaux as follows: given
$S,T\in \bT_n$ set  $T\ao S$ if $I_T\subset I_S.$ Respectively, we
define a geometric order on Young tableaux as follows: $S,T\in
\bT_n$ set  $T\go S$ if $\ov\Vscr_T\supset \ov\Vscr_S.$
\par
The most natural conjecture is that the geometric order on Young tableaux coincides
with the algebraic order.
The result about irreducibility of a variety associated
to a simple highest weight module provides the implication $T\ao S\ \Ra\ T\go S.$
Unfortunately we have no
algebro-geometrical tools to show the other implication.

\subsection{}\label{1.7}
Let us return to the description of an orbital variety
closure in a  semisimple $\gog.$
This description has two components.
The first purely geometrical component is what varieties constitutes the closure of an
orbital variety. This question can be formulated as following. Let $\Vscr$ be an orbital
variety then its $G$-saturation $\Oscr_{\Vscr}$  is  a nilpotent orbit, $\Vscr$ is associated
to. Let us take $\Oscr\st\ov\Oscr_{\Vscr}$ and consider $\ov\Vscr\cap \Oscr.$
As shown in \cite{M6}, this intersection is always not empty. Hence, a natural task
is to describe the irreducible components of this intersection. Is this intersection equidimentional?
Is this intersection Lagrangian?
\par
Again, as shown in \cite{M6}, if $\gog$ contains factors not of type $A_n$ there exist
orbital varieties in $\gog$ such that the intersection mentioned above
is not Lagrangian. However, the same argument does not work if all factors
are of type $A_n.$ As shown in Part I, 4.1.8, in that case $\ov\Vscr\cap \Oscr$ contains at
least one orbital variety. Moreover,
for some special cases  in $\gs\gl_n$ (cf. Part II, 2.3 and \cite[4.2]{M3})
the intersection is equidimentional and Lagrangian.
Together with the computations in low rank cases these facts support the conjecture
that in $\gs\gl_n$ the closure of an orbital variety is a union of orbital
varieties.

\subsection{}\label{1.8}
The other component of the description of an orbital variety closure
is combinatorial, that is the description of orbital varieties in the closure
of a given one in terms of Young tableaux only. We will discuss this in terms
of different partial orders.
Since we work with different partial orders and compare them
we will use the following terminology, customary in combinatorics.
Given two partial order relations on a set $S$
we call an order $\stackrel b {\leq}$ an extension of an order $\stackrel a {\leq}$ if
$x{\stackrel a {\leq}} y$ implies $x{\stackrel b {\leq}} y$ for any $x,y\in S.$
We will also call $\stackrel a {\leq}$ a restriction of  $\stackrel b {\leq}$ in that case. We
denote this by ${\stackrel a {\leq}}\preceq {\stackrel b {\leq}}.$
If $\stackrel b {\leq}$ is a proper extension of $\stackrel a {\leq}$
we write ${\stackrel a {\leq}}\prec {\stackrel b {\leq}}$

\subsection{}\label{1.9}
As we have already mentioned in Part I, the orbital varieties derive from
the works of N. Spaltenstein ~\cite{Sp1} and ~\cite{Sp2},
and R. Steinberg ~\cite{St1} and \cite{St2} during their studies of unipotent
variety of a complex semi-simple group $\bf G.$
\par
Recall the notion of
$L_w$ from \ref{1.5}. Any primitive ideal from $\Xscr_0$ is just
$I_w:=I(L_w),$ as it is explained in short in \ref{1.3}. Let us explain the results
of M. Duflo in more details.

Recall that each $w\in W$ is a product of fundamental reflections $s_{\al}\ :\ \al\in\Pi.$
We denote by $\ell (w)$ the minimal  length  of any such expression for $w.$
M. Duflo was the first (\cite{D}) to show that for any
semi-simple Lie algebra $\gog$ and its Weyl group $W$ if $w,y\in W$ are such that
$w=yx$ and $\ell(w)=\ell(x)+\ell(y)$
then $I_{w^{-1}}\supset I_{y^{-1}}.$  In that case we put $y\dor w$ and call it
a Duflo order. The more standard name for this order is a weak (right) Bruhat order.
However, because of the result, described above,
we prefer to call it a Duflo order in our context.
By Steinberg's construction
it was obvious that Duflo order implies the inclusion of orbital varieties as well,
that is, if $y\dor w$ then $\ov\Vscr_w\st \ov\Vscr_y$ (because of the inclusion of generating
subspaces). As we already mentioned in \ref{1.6}, the irreducibility
of an associated variety in $\gs\gl_n$ implies that if $I_{w^{-1}}\supset I_{y^{-1}}$ then
$\ov\Vscr_w\st \ov\Vscr_y$ so that in that case we do not need even Steinberg's construction
to show that induced Duflo order is the restriction of the geometric order.
\par
Induced Duflo order on Young tableaux was
the main object  of Part I. We denote it by $\dor.$
The purely combinatorial nature of the decomposition
into the cells, the above relation between the
geometric and the induced Duflo order and  the computations
for low rank cases
lead one to expect that  both the algebraic and the geometric orders
must coincide with the induced Duflo order.
However, this is false.
As we show
in \ref{5.6}, the induced Duflo order coincides with the algebraic and the geometric orders
for $n\leq 5$ and it is a proper restriction of
the algebraic  order (hence, also of the geometric order) for $n\geq 6.$
\par
Using Spaltenstein's construction we consider each orbital variety
as a double chain of nilpotent orbits (cf. \ref{3.9}). The inclusions on nilpotent orbit closures
are described combinatorially by Gerstenhaber's construction, explained in short in \ref{3.1}.
Thus, we can define another combinatorial order
on orbital varieties by inclusions of all nilpotent  orbit closures in the double chains.
We call it the chain order and denote by $\cho.$
This order was described in \cite{M0} and \cite{v-L}.
It is an extension of the geometric order. It coincides with the geometric order
for $n\leq 5$ and it is its proper extension for $n\geq 6.$
By a natural and very slight extension of the chain order one can force it to
coincide with the geometric order for $n=6,$ however, for $n\geq 7$  the new
chain order is a proper extension of the geometric order. We demonstrate this in \ref{3.10}.
Thus, for $n\geq 7$ the relations between the orders are
$${\dor\prec \ao \preceq \go\prec \cho}$$
We have to extend the induced Duflo order and to restrict the chain order to get two new orders
which will sandwich the algebraic and the geometric orders.

\subsection{}\label{1.10}
To do this we again return to the theory of primitive ideals of $U(\gs\gl_n).$
Here D. Vogan invented a beautiful
technique of an order isomorphism $\Tscr_{\al,\be}$ (cf. \cite{V2}).
Let us explain it in more detail for Young tableaux.
Recall that in the case of $\gs\gl_n$ one has $\Pi=\{\al_i\}_{i=1}^{n-1}.$
The notion of
$\tau(T)$ is defined as follows (cf. Part I, 2.4.14 for details).
For $T\in \bT_n$ and $a\ : 1\leq a\leq n$
we set $r_{\sr T}(a)$ to be the number
of a row $a$ belongs to. In these terms $\tau(T):=\{\al_i\ :\ r_{\sr T}(i+1)>r_{\sr T}(i)\}.$
Now let $\al,\be\in\Pi$ be subsequent roots (that is for $\al=\al_i$  $\be$ must be either $\al_{i+1}$
or $\al_{i-1}$). For such $\alpha,\beta$ put $D_{\al,\be}=\{T\in \bT_n\ :\ \al\not\in \tau(T),\ \be\in\tau(T)\}.$
Vogan's bijection $\Tscr_{\al,\be}$ maps $D_{\al,\be}$ onto $D_{\be,\al}.$
For  $T\in D_{\al,\be}$ we obtain $\Tscr_{\al,\be}(T)$ by changing
numbers in two boxes. We explain this purely combinatorial
procedure in \ref{5.2}. By \cite{V2} this  procedure preserves the algebraic order, that is
for $T,S\in D_{\al,\be}$ one has $T\ao S$ if and only if
$\Tscr_{\al,\be}(T)\ao \Tscr_{\al,\be}(S).$ A. Joseph showed in \cite{Jo1} that $\Tscr_{\al,\be}$
can be applied to orbital varieties as well. Slightly generalizing his result, we get in \ref{5.5}
that $\Tscr_{\al,\be}$ preserves the geometric order as well.
\par
Further we show that both the induced Duflo and the chain orders are not preserved under $\Tscr_{\al,\be}.$
These facts provide us examples showing that the induced Duflo order is a proper restriction and
the chain order is a proper extension of both the algebraic and the geometric orders.

\subsection{}\label{1.11} Moreover, we use $\Tscr_{\al,\be}$ to extend the induced Duflo
order and to restrict the chain order as we explain in short in this section.

As it was shown  by A. Joseph (cf. \cite[3.12]{B-V}),
Robinson - Schensted insertion (cf. Part I, 2.4.5, 2.4.10) preserves the algebraic order, that is
$T\ao S$ implies both $(T\Da a)\ao (S\Da a)$ and $(a\Ra T)\ao (a\Ra S).$  We show in
\ref{6.6}  that the same is true for the geometric order.
This gives us in particular a geometric meaning of
Robinson-Schensted procedure. This we believe explains why
Robinson-Schensted procedure describes the cell decomposition.
\par
Moreover, the above property of Robinson-Schensted insertion
together with $\Tscr_{\al,\be}$ operators
leads to an extension of induced Duflo order defined in \ref{6.7}
We call it Duflo-Vogan order and denote it by $\dvo$. It is a restriction of the algebraic order.
\par
On the other hand,  we use $\Tscr_{\al,\be}$ to restrict the chain order (cf. \ref{6.7}).
We call this restriction Vogan-chain order. It is an extension of the  geometric order.
\par
Now we have two combinatorially defined orders (however, of very different nature)
and both the algebraic and the
geometric orders are sandwiched between them. Computer computations show that
they coincide for $n\leq 9.$ In the case  of $n=10$ there is one example
(up to operations $\Tscr_{\al,\be}$ and transposition) of $T,S\in \bT_{10}$ such that
$T\chvos S$ however, $T{\stackrel{\rm D-V}{\not<}} S.$ In this case we check the inclusion of
primitive ideals
with the help Kazhdan-Lusztig combinatorics (using, in particular, the program Coxeter 3.0 $\beta 2$
of F. Du Cloix for the computations of Kazhdan-Lusztig polynomials).
The computations show that
$T\aos S$. Thus, in that case we have that the algebraic the geometric and Vogan-chain orders coincide.
All these facts support our main conjecture.
\begin{conj}
For any $T,S\in \bT_n$ $\ov\Vscr_T\st\ov\Vscr_S$ if and only if $I_T\supset I_S$
and this happens if and only if $S\chvo T.$
\end{conj}

\subsection{}\label{1.12} Given a set $S$ and a partial order $\leq$ on it, the cover of $a\in S$
for this order
is a set of all $b\in S$ such that $b>a$ and for any $c\in S$ such that $b\geq c\geq a$
one has either $c=b$ or $c=a.$ To describe the cover of an element for a given order is a delicate
question even in the cases when we have a satisfactory description of the order. In Part I we discussed
the cover of a tableau
for the induced Duflo order. As we have shown in Part II all our orders coincide for Richardson component.
The full description of the cover of a Richardson orbital variety is provided in Part II. Here we
discuss the cover of a tableau for
the geometric order. However, here our results are mostly of negative nature. The only positive result is
that $\Tscr_{\al,\be}$ preserves the cover for the geometric and algebraic orders, i.e., let
$T,S\in D_{\al,\be}$ then $S$ is in the cover $T$
iff $\Tscr_{\al,\be}(S)$ is in the cover of $\Tscr_{\al,\be}(T).$ On the other hand
neither projection, nor injection preserve
the cover. As well, we show that $\Vscr_S$ being in the cover of $\Vscr_T$ does not imply that
the nilpotent orbit $\Oscr_{\Vscr_S}$ is in the cover
of $\Oscr_{\Vscr_T}.$ This again demonstrates that the description of  inclusion of orbital variety closures
is a much more delicate problem than the description of  inclusion of nilpotent orbits.

\subsection{}\label{1.13}
The body of the paper consists of  6 sections.
\par
Sections 2-4 are preliminary. For the convenience of the reader we repeat necessary notation and
results from Parts I and II
which can be formulated in short. If the formulation is too long, as for example, in the
case of Robinson-Schensted procedure we provide the exact reference to the subsection
of Parts I and II. I hope these  sections make the paper self-contained.

Section 2 is devoted to the facts connected to Steinberg's construction and the induced Duflo order
essential in further analysis. In section 3 we explain
Spaltenstein's construction,  the connection between
Steinberg's and Spaltenstein's construction for $\gs\gl_n$ and define the chain order.
In section 4 we consider the facts from the theory of primitive ideals
essential in  the subsequent analysis and consider the algebraic order.

In section 5 we consider Vogan's $\Tscr_{\al,\be}$ operator and show that it is a geometric
order isomorphism. Section 6 is devoted to the description of
an orbital variety closure and the comparison of different orders on Young tableaux.
Finally, in section 7 we discuss the questions connected to the cover of a given tableau for
the geometric order.
\par
In the end one can find the index of notation in which
symbols appearing frequently are given with the subsection
where they are defined. We hope that this will help the reader
to find his way through the paper.
\parno
{\bf Acknowledgments.}\ \ \ I would like to express my gratitude to A. Joseph for
introducing to me the world of orbital varieties, for proposing the
ideas underlying the chain order and many fruitful discussions through the various
stages of this work.

I would like to thank V. Hinich
for many fruitful discussions on algebraic geometry connected to this research,
and F. Du Cloix for providing me with his package of programs
``Coxeter 3.0 $\beta 2.$ Without this package the computations of Kazhdan-Lusztig polynomials
in $\gs\gl_{10}$ would be impossible.

\section{Steinberg's construction and induced Duflo order}

\subsection{}\label{2.1}
In this section we repeat the definitions and facts from Part I that
we need in our further discussion.

Let us explain in short Steinberg's construction of orbital varieties.
In detail it is described in Part I, \S 2.1.2, 2.1.3.

Let $\gog$ be any semisimple Lie algebra. Fix its triangular
decomposition $\gog=\gn\oplus\gh\oplus\gn^-$. Let $\bB$ be the Borel
subgroup of $\bG$ with $\Lie(\bB)=\gh\oplus \gn$ and let $\bB$ act
adjointly on $\gn.$ Recall notation for root system from \ref{1.5}.
Let $X_\alpha$ denote the root subspace for $\alpha\in R.$ One has
$\gn=\bigoplus\limits_{\alpha\in R^+}X_\alpha.$ Let $W$ be the Weyl
group of $<\gog,\gh>.$ The action of $w\in W$ on root subspace
$X_\alpha$ is defined (in a standard way) by
$w(X_\alpha)=X_{w(\alpha)}.$

For $H\subset \bG$ and $\ga\subset\gog$ put $H(\ga):=\{AxA^{-1}\ :\
A\in H,\ x\in\ga\}$. Let $\ov\ga$ mean the closure of $\ga$ (in
Zarisky topology).

Consider the following subspace of $\gn:$
$$\gn\cap^w\gn=\bigoplus\limits_{\{\alpha\in R^+\, |\, w^{-1}(\alpha)\in R^+\}}X_\alpha.$$
 Consider $\ov{\bG(\gn\cap^w\gn)}.$ Since the
number of orbits is finite this is a closure of the unique orbit
which we denote by $\Oscr_w.$ By Steinberg \cite{St1} one has
\begin{theorem} For each $w\in W$ there exists an orbital variety
$\Vscr$ and  for each orbital variety
$\Vscr$ there exists $w\in W$ such that
$$\Vscr=\ov{\bB (\gn\cap{^w}\gn)}\cap\Oscr_w.$$
\end{theorem}

In what follows we will denote $\Vscr_w:=\Vscr$ in that case.

\subsection{}\label{2.2}
As we have explained in \ref{1.2} the Weyl group $W$ is partitioned into
geometric cells according to Steinberg's construction:
for $w\in W$ let $\Cscr_w=\{y\in W\ |\ \Vscr_y=\Vscr_w\}.$

To give a description of geometric cells in $\gs\gl_n$ we need the notion
of Young tableaux and Robinson-Schensted procedure.

Recall that the Weyl group of $\gs\gl_n$ is $\bS_n.$ We will write
elements of $w\in W$ in a word form, that is $w=[a_1,\ldots,a_n]$
means $w(i)=a_i$ for $1\leq i\leq n.$ In this case
$\Pi=\{\al_i\}_{i=1}^{n-1}$ and transpositions
$s_i:=s_{\al_i}=[1,\ldots,i+1,i,\ldots,n]$ are generators of $\bS_n$
as a Weyl group.

Given a  partition $\lambda$ of $n$ written in decreasing order
$\lambda=(\lambda_1\geq \lambda_2\geq\ldots \geq\lambda_m>0)$
we define Young digram $D_{\lambda}$ to be an array of $m$  rows of boxes
starting on the left  with the $i$-th row containing $\lambda_i$ boxes.
Set $\bD_n$ to be the set of all Young diagrams corresponding to $n.$ (By a slight abuse
of notation we will not distinguish in what follows between a Young diagram $D_{\lambda}$
the partition $\lambda.$)

Given a Young diagram $D_{\lambda}\in \bD_n$ we can fill its boxes with the integers
$1,\ldots,n.$ If numbers increase in rows from left to right and in columns from top to bottom,
such an array is called a (standard) Young tableau of shape $\lambda.$ Given a Young tableau $T$ let
$\sh(T)$ denote the corresponding  partition.
Set $\bT_n$ to be the set of all Young tableaux
with $n$ boxes and put $\bT_{\lambda}$ to be the set of all Young tableaux of shape $\lambda.$

Given $w\in\bS_n$ let $(P(w),Q(w))$
be the pair of Young tableaux (of the same shape) constructed with the help of
Robinson-Schensted procedure.
(cf., for example \cite[\S 3]{Sa}).

On the other hand, partitions are connected to the nilpotent orbits
in a natural way via Jordan form. Recall that $\bG=SL_n$ acts on
$\gog$ by conjugation. For any $u\in\gog$ its $\bG-$orbit $\Oscr_u$
is defined by Jordan form $J(u)$. Let $\Nscr$ denote the nilpotent
cone in $\gog.$ If $u\in \Nscr$ all its eigenvalues are 0 and Jordan
form of $u$ is defined only by the length of its Jordan blocks.
Writing the lengths of Jordan blocks in decreasing order we get a
natural bijective correspondence between nilpotent orbits of
$\gs\gl_n$ and partitions of $n$ (also written in decreasing order).
The map $u\rar J(u)$ gives a bijection of $\Nscr/\bG$ onto $\bD_n.$
Given $J(u)=\lambda$  we also write $\Oscr_{\lambda}:=\Oscr_u$ and
$D_u:=D_{\lambda}.$

For $\gog=\gs\gl_n$ as it is shown in \cite{St2} one has
\begin{theorem}
\begin{itemize}
\item[(i)] $\Cscr_w=\{y\in \bS_n\ |\ P(y)=P(w)\};$
\item[(ii)] $\Vscr_w$ is associated to $\Oscr_{\sh(P(w))}.$
\end{itemize}
\end{theorem}

\subsection{}\label{2.3}
Let us return to Duflo order described in \ref{1.9}. By definition
one can see immediately that $y\dor w$ implies $\gn\cap^y\gn\supset \gn\cap^w\gn.$
Therefore $y\dor w$ implies by Steinberg's construction $y\go w.$
We induce this order to the order on the set of orbital varieties and respectively
to the set of Young tableaux. We continue to denote it by $\dor.$

As we show in \ref{5.6}, for $n\geq 6$ induced Duflo order is a proper restriction of
both the algebraic and the geometric orders.

\subsection{}\label{2.4}
In what follows we will need some theorems from Part I, \S 4.1.1 on projections of orbital variety
closures onto Levi factor. To formulate them we recall the definitions and notation from  Part I,
\S 2.1.8.

For any $\alpha\in \Pi$ let $\bP_\alpha$ be the standard parabolic
subgroup of $\bG$ such that $Lie(\bP)=\gn\oplus\gh\oplus
X_{-\alpha}.$

Given $\Iscr\st \Pi,$ let ${\bP}_{\Iscr}$ denote the unique
standard parabolic subgroup of $\bG$ such that
$\bP_{-\alpha}\st \bP$ iff $\alpha\in\Iscr.$ Let ${\bM}_{\Iscr}$ be the unipotent
radical of ${\bP}_{\Iscr}$ and ${\bL}_{\Iscr}$ a Levi factor.
Let $\gp_{\sr \Iscr},\ \gm_{\sr \Iscr},\ \gl_{\sr \Iscr}$ denote
the corresponding Lie algebras.
Set ${\bB}_{\Iscr}:={\bB}\cap{\bL}_{\Iscr}$ and
$\gn_{\sr \Iscr}:=\gn\cap\gl_{\sr \Iscr}.$
 We have decompositions ${\bB}={\bM}_{\Iscr}\ltimes{\bB}_\Iscr$ and
$\gn=\gn_{\sr \Iscr}\oplus\gm_{\sr \Iscr}.$ They define projections
${\bB}\rar {\bB}_\Iscr$ and $\gn\rar\gn_{\sr \Iscr}$ which we denote by
$\pi_{\sr \Iscr}.$

Set $W_\Iscr:=<s_{\al}\ :\ \al\in \Iscr>$ to be a
parabolic subgroup of $W.$ Set
$F_\Iscr:=\{w\in W\ :\ w(\al)\in R^+\ \forall \ \al\in \Iscr\}$ and
$F^{-1}_\Iscr:=\{w^{\sr -1}\ :\ w\in F_\Iscr\}.$
Set $R_\Iscr^+=R^+\cap \Spa(\Iscr).$
A well-known classic result (cf., for example \cite{Ca}) is that
each $w\in W$ has a unique expression of the
form $w=w_{\sr \Iscr}f_{\sr \Iscr}$
where $f_{\sr \Iscr}\in F^{-1}_\Iscr,\ w_{\sr \Iscr}\in W_\Iscr$ and $\ell(w)=\ell(w_{\sr \Iscr})+
\ell(f_{\sr \Iscr}).$  Moreover,
$$R_\Iscr^+\cap^w R^+=R^+_\Iscr\cap^{w_{\sr \Iscr}} R^+_\Iscr.$$
Thus, decomposition $W=W_\Iscr\times F^{-1}_\Iscr$
defines a projection $\pi_{\sr \Iscr}:W\rar W_\Iscr.$ For $w\in W$ set $w_{\sr \Iscr}:=\pi_{\sr \Iscr}(w).$
This element can be regarded as an element of $W_\Iscr$ and as an element of $W.$

Let ${\Cscr}_{w_{\sr \Iscr}}$ denote its cell in $W$ and ${\Cscr}^{\Iscr}_{w_{\sr \Iscr}}$ denote
its cell in $W_\Iscr.$ Respectively let ${\Vscr}_{w_{\sr \Iscr}}$ be the corresponding orbital variety
in $\gog$ and ${\Vscr}^{\Iscr}_{w_{\sr \Iscr}}$ be the corresponding orbital variety
in $\gl_{\sr \Iscr}.$ All the projections
are in correspondence on orbital varieties and cells, namely
\begin{theorem} Let $\gog$ be a reductive algebra. Let $\Iscr\st \Pi.$
\begin{itemize}
\item[(i)]  For every $w\in W$ one has  $\pi_{\sr \Iscr}(\Cscr(w))\st \pi_{\sr \Iscr}(\Cscr(w_{\sr \Iscr}))=
            {\Cscr}_{\Iscr}(w_{\sr \Iscr}).$
\item[(ii)] For every orbital variety ${\Vscr}_w\st \gog$ one has $\pi_{\sr \Iscr}(\ov{\Vscr}_w)=
\ov{{\Vscr}^{\Iscr}_{w_\Iscr}}.$
\end{itemize}
\end{theorem}
By theorem we get that $w\go y$  implies that $w_{\sr \Iscr}\go y_{\sr \Iscr}$ both as elements of $W$
and as elements of $W_\Iscr$ for any $\Iscr\st\Pi.$

\subsection{}\label{2.5}
Let us list a few elementary properties of induced Duflo order. They are true in general
but we formulate them only for $\gs\gl_n$ since they are expressed nicely in terms of
Young tableaux.

We begin with a  well known result, shown for example in \cite[2.3]{JM}
\begin{prop} For any $w,y\in \bS_n$ one has $w\dor y$ iff for any
$i,j\ :\  1\leq i<j\leq n$ if
$w^{-1}(i)>w^{-1}(j)$ then  $y^{-1}(i)>y^{-1}(j).$
\end{prop}

\subsection{}\label{2.6}
We use here a few classical algorithms on Young tableaux which we
describe below. But first we need to set up the notation. Given
words $w=[a_{\sr 1},\ldots,a_i]$ and $y=[b_{\sr 1},\ldots, b_j]$
such that $\{a_s\}_{s=1}^i\cap\{b_s\}_{s=1}^j$ put $[w,y]=[a_{\sr
1},\ldots,a_i, b_{\sr 1},\ldots, b_j]$ to be their colligation.
Given a word $w=[a_1,\ldots, a_n]\in \bS_n$ and $a\in
\{1,\ldots,n+1\}$ put $\ov w_a:=[b_{\sr 1},\ldots,b_n]$ where
$$b_i=\begin{cases}
a_i& {\rm if}\  a_i<a,\\
a_i+1& {\rm otherwise.}\\
\end{cases}
$$
Note that by proposition \ref{2.5} we get immediately that $w\dor y$
provides $[\ov w_a,a]\dor[\ov y_a,a]$ and $[a,\ov w_a]\dor[a,\ov
y_a].$

 Let
$T$ be a Young tableau. We denote the content of a box on the
intersection of $i-$th row and $j-$th column by $(T)_{i,j}.$

Given a tableau $T$ (such that its elements are some integers among
$1,\ldots,n$ but not all of them) and an integer $a\in
\{i\}_{i=1}^n$ which is  not among the elements of $T,$ then the
Robinson-Schensted insertions $(T\Da a),\ (a\Ra T)$ can be defined
(cf. Part I, 2.4.5, 2.4.10, or \cite[\S 3]{Sa}). This procedure
gives us (inductively) $w\rar P(w).$

Given a Young tableau $T$ put and $a\in \{1,\ldots,n+1\}$ let $\ov
T_a$ be obtained from $T$ by
$$(\ov T_a)_{i,j}=\begin{cases}
(T)_{i,j}& {\rm if}\  (T)_{i,j}<a,\\
(T)_{i,j}+1& {\rm otherwise.}\\
\end{cases}
$$
By Robinson-Schensted procedure one has that $P([\ov
w_a,a])=(\ov{P(w)}_a \Da a)$ and $P([a,\ov w_a])=(a\Ra
\ov{P(w)}_a).$

Thus, by the previous note if $S,T\in \bT_n$ are such that $S\dor T$
then for any $a\in\{i\}_{i=1}^{n+1}$ one has $(a\Ra \ov S_a)\dor
(a\Ra \ov T_a)$ and $(\ov S_a\Da a)\dor (\ov T_a\Da a).$

For any $i,j\ :\ 1\leq i<j\leq n$ set
$\Iscr_{i,j}=\{\alpha_k\}_{k=i}^{j-1}$ and set $\pi_{i,j}:=\pi_{\sr
\Iscr_{i,j}}.$ In that case for $w=[a_1,\ldots,a_n]$ one can
consider $\pi_{i,j}(w)$ as a word in $\bS_{i,j}$ -- symmetric group
of $i,\ldots,j$ then $\pi_{i,j}(w)$ is obtained simply be deleting
$1,\ldots, i-1$ and $j+1,\ldots,n$ in the word $[a_1,\ldots,a_n].$

Let us recall that given a Young diagram
For $T\in\bT_n$ let $\pi_{i,j}(T)$ be the new tableau
with entries $i,\ldots,j$ obtained from $T$ by   Sch\"utzenberger's ``jeu de taquin''
process. Then by Sch\"utzenberger one has
$$P(\pi_{i,j}(w))=\pi_{i,j}(P(w)).\eqno{(*)}$$
All the details
can be found in \cite[\S 3]{Sa}.
As a straightforward corollary of proposition \ref{2.5} and construction of $\pi_{i,j}$
we get that  inequality $w \dor y$ implies $\pi_{i,j}(w)\dor
\pi_{i,j}(y)$ for any $i,j\ :\ 1\leq i<j\leq n.$

Given a Young tableau $T$ let $T^{\dagger}$ denote the
transposed tableau, that is a tableau obtained from the given one by interchanging rows and columns.
For example
$$
\vcenter{
\halign{& \hfill#\hfill
\tabskip4pt\cr
\multispan{7}{\hrulefill}\cr
\ssa
\vb& 1 & &3  & & 5 &\ts\vb\cr
\vsa
&&&&\multispan{3}{\hrulefill}\cr
\ssa
\vb& 2 & &4 &\ts\vb\cr
\vsa
\multispan{5}{\hrulefill}\cr}}^{\dagger}=
\vcenter{
\halign{& \hfill#\hfill
\tabskip4pt\cr
\multispan{5}{\hrulefill}\cr
\ssa
\vb& 1 & &2  & \ts\vb\cr
\vsa
&&&&\cr
\ssa
\vb& 3 & &4 &\ts\vb\cr
\vsa
&&\multispan{3}{\hrulefill}\cr
\ssa
\vb& 5 &\ts\vb\cr
\vsa
\multispan{3}{\hrulefill}\cr}}
$$

Note that $w_o=[n,n-1,\ldots,1]$ is the maximal element of $\bS_n$ in Duflo order.
Obviously, for any
$w=[a_{\sr 1}, a_{\sr 2},\ldots,a_n]$ one has $ww_o=[a_{n}, a_{n-\sr 1},\ldots,a_{\sr 1}].$
Moreover, by the results of Sch\"utzenberger,
$P(ww_o)=(P(w))^{\dagger}$ (cf. \cite[\S 3]{Sa}).
As a straightforward corollary of proposition \ref{2.5} and this construction we get
that $w\dor y$ iff $ww_o\dg yw_o.$

Summarizing this subsection in terms of Young tableaux we get
\begin{cor} For any Young tableaux $S,T\in \bT_n$ one has
\begin{itemize}
\item[(i)] If $S\dor T$ then for any $i,j\ : 1\leq i<j\leq n$ $\pi_{i,j}(S)\dor\pi_{i,j}(T);$
\item[(ii)] If $S\dor T$ then for any $a\in \{i\}_{i=1}^{n+1}$ one has
$(\ov S_a\Downarrow a)\dor (\ov T_a\Downarrow a)$ and $(a\Ra\ov
S_a)\dor (a\Ra \ov T_a).$
\item[(iii)] $S\dor T$ iff $S^{\dagger}\dg T^{\dagger}.$
\end{itemize}
\end{cor}

Note that for $w=[a_{\sr 1}, a_{\sr 2},\ldots,a_n]$ one has $w_ow=[n+1-a_{\sr 1},\ldots,n+1-a_n]$
so that again by proposition \ref{2.5} $w\dor y$ iff $w_ow\dg w_oy.$ However,
in this case $P(w_ow)$ is expressed in a more complex way using ``evacuation'' procedure of
Sch\"utzenberger (cf. for example \cite[\S 3]{Sa}). We are not going to discuss this procedure here.

\section{Spaltenstein's construction and chain order}

\subsection{}\label{3.1}
We begin with the construction of
Gerstenhaber giving the combinatorial description of
nilpotent orbit closure (in $\gs\gl_n$). It is described in detail
in many places including Part I \S 2.3.

Recall from \ref{2.2} that the orbits of elements of nilpotent
cone $\Nscr$ under the action of conjugation by $SL_n$ are completely
described by Young diagrams via the Jordan form. Recall the notation
$\Oscr_{\lambda}$ from \ref{2.2}.

We define  an order relation on Young diagrams as follows.
Let $D_{\lambda}=(\lambda_{\sr 1},\cdots,\lambda_k)$ and $D_{\mu} = (\mu_{\sr 1},\cdots,\mu_j) $
be Young  diagrams in  $\bD_n$. Set $D_\lambda \geq D_\mu$ if for each
$i\ : 1\leq i\leq \min (k,j)$ one has
$$ \sum^i_{m=1}\lambda_m \leq \sum^i_{m=1}\mu_m \ . $$

(Usually the order relation goes the other way round, but to put it
in correspondence with the inclusions on primitive ideals we choose
this direction.)

Then by  Gerstenhaber (cf. \cite[3.10]{H}, for example) one has

\begin{theorem} Let $\mu$ be a partition of $n$
and $\Oscr_\mu$ be the corresponding nilpotent orbit in
$\gs\gl_n.$ Then
$$ \overline{\Oscr_\mu} = \coprod_{\lambda\vert D_\lambda\geq D_\mu}
                                                        \Oscr_\lambda \ . $$
\end{theorem}

\subsection{}\label{3.4}
Finally we explain Spaltenstein's construction (\cite{Sp1}).
Let $\gog=\gs\gl_n$ and let $\gn=\gn_n$ be the subspace of strictly upper triangular
matrices. For any $u\in \gn$ let $J(u)$ denote the partition corresponding to
its Jordan form as in \ref{3.1}.
One can consider the projection $\pi_{1,k}:\gn_n \rar \gn_k$ obtained by deleting
rows and columns $k+1,\ldots,n.$ For $k\ :\ 1\leq k\leq n-1$  put $u_{1,k}=\pi_{1,k}(u)$
and put $u_{1,n}=u.$ For any $u$ we construct the chain of Young diagrams by
$\theta(u)=\{J(u_{1,n}),J(u_{1,n-1}),\ldots,J(u_{1,1})\}.$
For example,
$$u=\left(\begin{array}{cccccc}
0 & 0 & 0 & 0 & 1 & 0 \\
0 & 0 & 1 & 0 & 0 & 0 \\
0 & 0 & 0 & 0 & 0 & 1 \\
0 & 0 & 0 & 0 & 0 & 0 \\
0 & 0 & 0 & 0 & 0 & 1 \\
0 & 0 & 0 & 0 & 0 & 0 \\
\end{array}\right),\qquad \theta(u)=\{ (3,2,1),\, (2,2,1),\, (2,1,1),\, (2,1),\, (1,1),\, (1)\}.$$

Recall that given a Young diagram $D_\lambda$ where $\lambda$ is a partition of $n$ we
construct a (standard) Young tableau (associated to $D_\lambda$) by filling in all the boxes with numbers $1,\ldots,n$
in such a way that entries increase in rows from left to right and in columns from top to bottom.
If $T$ is a Young tableau associated to $D_\lambda$ we will denote $\sh(T)=\lambda.$

 Define a projection $\pi_{1,k}:\bT_n\rar \bT_k$ by removing cells
containing numbers $k,k+1,\ldots,n.$ Note that $\pi_{1,k-1}(T)$ is
obtained from $\pi_{1,k}(T)$ by deleting exactly one box containing
$k.$ So if we know $\pi_{1,k-1}(T)$ and $\sh(\pi_{1,k}(T))$ we can
reconstruct $\pi_{1,k}(T)$ by putting $k$ in the only new box of
$\sh(\pi_{1,k}(T)).$ In such a manner we get a bijection between
$\bT_n$  and the set of chains of Young diagrams such that each
diagram in the chain differs from the previous diagram by one box.
Set $\phi(T)=\{\sh(T),\,
\sh(\pi_{1,n-1}(T)),\ldots,\sh(\pi_{1,1}(T))\}.$ For example
$$T=\vcenter{
\halign{& \hfill#\hfill
\tabskip4pt\cr
\multispan{7}{\hrulefill}\cr
\ssa
\vb& 1 & &3  & & 6 &\ts\vb\cr
\vsa
&&&&\multispan{3}{\hrulefill}\cr
\ssa
\vb& 2 & &5 & \ts\vb\cr
\vsa
&&\multispan{3}{\hrulefill}\cr
\ssa
\vb&4 & \ts\vb\cr
\vsa
\multispan{3}{\hrulefill}\cr}}\quad \phi(T)= \{(3,2,1),\, (2,2,1),\, (2,1,1),\, (2,1),\,
(1,1),\, (1)\}.$$
For $T\in \bT_n$ set $\nu_{\sr T}:=\theta^{-1}_T=\{u\in\gn\ |\ \theta(u)=\phi(T)\}.$
In our examples of this subsection we have $u\in\nu_{\sr T}.$

Given $\lambda \vdash n$ let $\Oscr_\lambda$ be the corresponding nilpotent
orbit in $\gs\gl_n$ and let $\bT_\lambda$ be the set of standard Young tableaux of shape
$\lambda$ (as it is defined in \ref{2.2}). Then as a straightforward corollary of
the construction we get $\Oscr_\lambda\cap\gn=
\coprod\limits_{T\in\bT_\lambda}\nu_{\sr T}.$ As it was shown in \cite{Sp1}

\begin{theorem} For any orbital variety $\Vscr$ associated to $\Oscr_\lambda$ there exists
$T\in\bT_\lambda$ such that $\nu_{\sr T}$ is a dense open part of $\Vscr$
and for any
$T\in\bT_\lambda$ there exists an orbital variety $\Vscr$ associated to $\Oscr_\lambda$
such that $\ov{\nu_{\sr T}}\cap\Oscr_\lambda=\Vscr.$
\end{theorem}
In what follows we will denote $\Vscr_T:=\Vscr$ in that case.

\subsection{}\label{3.7}
At first we will use theorem \ref{2.4} to show a well known fact
that Spaltenstein's and Steinberg's constructions give exactly the same orbital varieties.
\begin{prop}
Given $T\in \bT_n,$ let $w_{\sr T}\in \bS_n$ be some element such that $P(w_{\sr T})=T.$
Let $\Vscr_{w_{\sr T}}$ and $\Vscr_T$ be the corresponding orbital varieties.
Then $\Vscr_{w_{\sr T}}=\Vscr_T.$
\end{prop}
\Pf
Indeed, the claim is trivially true for $\gs\gl_2.$ Assume it is true for
$\gs\gl_{n-1}$ and show it for $\gs\gl_n.$

Let $T\in\bT_n$ be some tableau of shape $\lambda.$ Let $w\in\bS_n$
be such that $P(w)=T.$ Then by Theorem \ref{2.1} $\Vscr_w$ is some
orbital variety associated to $\Oscr_\lambda.$ Consider also
$\Vscr_T$ which is also associated to $\Oscr_\lambda$ just by
Spaltenstein's construction.

By \ref{2.4} and \ref{2.6} one has that
$\pi_{1,n-1}(\ov\Vscr_w)=\ov\Vscr_{\pi_{1,n-1}(w)}.$ Note also that
$P(\pi_{1,n-1}(w))=\pi_{1,n-1}(P(w))$ by \ref{2.6} $(*).$ Thus,
$\ov\Vscr_{\pi_{1,n-1}(w)}=\ov\Vscr_{\pi_{1,n-1}(T)}$ by induction
hypothesis.

Assume $\Vscr_w=\Vscr_{T\pr}$. Then by Theorem \ref{2.2}
$\sh(T\pr)=\sh(T)$ and by induction hypothesis
$\pi_{1,n-1}(T\pr)=\pi_{1,n-1}(T).$ Thus, $\phi(T)=\phi(T\pr)$ so
that $T\pr=T.$
 \QED

\subsection{}\label{3.8}
Moreover, in the same way we can get a more refine result.
\begin{prop}
Given $T\in \bT_n,$ let $w_{\sr T}\in \bS_n$ be some element such that $P(w_{\sr T})=T.$
Then $(\gn\cap^{w_T}\gn)\cap\nu_{\sr T}$ is dense in $\gn\cap^{w_T}\gn$ and
$\bB(\gn\cap^{w_T}\gn)\cap\nu_{\sr T}$ is dense in $\nu_{\sr T}.$
\end{prop}
\Pf
Indeed, since $\nu_{\sr T}$ is locally closed (in Zariski topology) and $\bB$ stable one has that
$(\gn\cap^{w_T}\gn)\cap \nu_{\sr T}\ne\emptyset$ (otherwise
$\bB(\gn\cap^{w_T}\gn)\st \ov\Vscr_T\setminus\nu_{\sr T}$ and thus, $\ov{\bB(\gn\cap^{w_T}\gn)}\st
\ov\Vscr_T\setminus\nu_{\sr T}$ which contradicts proposition \ref{3.7}.
Now since $\gn\cap^{w_T}\gn$ is closed in Zariski topology this implies that
$(\gn\cap^{w_T}\gn)\cap\nu_{\sr T}$ is dense in $\gn\cap^{w_T}\gn.$ Therefore $\ov
{\bB(\gn\cap^{w_T}\gn)\cap\nu_{\sr T}}=\ov\nu_{\sr T}$ which gives us that
$\bB(\gn\cap^{w_T}\gn)\cap\nu_{\sr T}$ is dense in $\nu_{\sr T}.$
\QED

\subsection{}\label{3.9}
Proposition \ref{3.8} together with theorem \ref{2.4} gives an idea of a  generalization of Spaltenstein's
construction.

Let $\hat\nu_{\sr T}:=\{u\in\gn\ :\
J(\pi_{i,j}(u))=\sh(\pi_{i,j}(T))\ \forall\ 1\leq i<j\leq n\}$. If
$u\in \gn\cap^{w_T}\gn$ is a generic matrix then  $\pi_{i,j}(u)$ is
a generic matrix of
$\pi_{i,j}(\gn)\cap^{w_{\pi_{i,j}(T)}}\pi_{i,j}(\gn)$ for any $i,j\
:\ 1\leq i<j\leq n$ so that $\pi_{i,j}(u)\in
\Oscr_{\sh(\pi_{i,j}(T))}$. Thus,
$(\gn\cap^{w_T}\gn)\cap\hat\nu_{\sr T}\ne\emptyset.$ Then exactly in
the same way as in \ref{3.8} we get that
$(\gn\cap^{w_T}\gn)\cap\hat\nu_{\sr T}$ is dense in
$\gn\cap^{w_T}\gn.$ As well one has $\hat\nu_{\sr T}\st \nu_{\sr T}$
and $\bB$ stable. Therefore we get $\ov{\hat\nu_{\sr
T}}=\ov\Vscr_T.$

In such a way we, generalizing Spaltenstein's construction, consider each Young tableau as a
double chain of Young diagrams. Put
$$\varphi(T):=\begin{array}{lcccc}
\sh(T)&          \sh(\pi_{1,n-1}(T))&\ldots&\sh(\pi_{1,2}(T))&\sh(\pi_{1,1}(T))\\
\sh(\pi_{2,n}(T)&\sh(\pi_{2,n-1}(T))&\ldots&\sh(\pi_{2,2}(T))&\\
\vdots&\ddots&&&\\
\sh(\pi_{n,n}(T))&&&&\\
\end{array}
$$
For example
$$\varphi\left(
\vcenter{
\halign{& \hfill#\hfill
\tabskip4pt\cr
\multispan{7}{\hrulefill}\cr
\ssa
\vb& 1 & &3  & & 4 &\ts\vb\cr
\vsa
&&&&\multispan{3}{\hrulefill}\cr
\ssa
\vb& 2 & &6 & \ts\vb\cr
\vsa
&&\multispan{3}{\hrulefill}\cr
\ssa
\vb&5 & \ts\vb\cr
\vsa
\multispan{3}{\hrulefill}\cr}}\right)=
\begin{array}{llllll}
(3,2,1)&(3,1,1)&(3,1)&(2,1)&(1,1)&(1)\\
(3,2)&(3,1)&(3)&(2)&(1)&\\
(2,2)&(2,1)&(2)&(1)&&\\
(2,1)&(1,1)&(1)&&&\\
(2)&(1)&&&&\\
(1)&&&&&\\
\end{array}
$$
\subsection{}\label{3.10}
Note that by theorem \ref{2.4}  $\ov\Vscr_1\st\ov\Vscr_2$ implies
$\pi_{\sr \Iscr}(\ov\Vscr_1)\st\pi_{\sr \Iscr}(\ov\Vscr_2)$ for any
$\Iscr\st\Pi.$\ This in turn implies the inclusion of corresponding
nilpotent orbit closures in $\gl_{\sr \Iscr}.$

Let $\Pi=\{\al_i\}_{i=1}^n$ be the set of simple roots in some simple Lie algebra.
For any connected $\Iscr\st \Pi$
let $\pi_{\sr \Iscr}$ be the corresponding projection. Set $\Oscr_{\pi_{\sr \Iscr}(w)}$ to be the nilpotent
orbit of $\pi_{\sr \Iscr}(w)$ in $\gl_{\sr \Iscr}.$

We define a partial order on orbital varieties and on $W$ as following
\begin{defi} Let $\gog$ be some simple Lie algebra.
For $y,w\in W$ (resp. for orbital varieties $\Vscr_y, \Vscr_w$)
set $y\cho w$ (resp. $\Vscr_y\cho \Vscr_w$) if
\begin{itemize}
\item[(i)] for any connected $\Iscr\st \Pi$
one has $\Oscr_{\pi_{\Iscr}(w)}\st \ov\Oscr_{\pi_{\Iscr}(y)};$
\item[(ii)] if for some $\Iscr$ $\ov\Oscr_{\pi_{\Iscr}(w)}= \ov\Oscr_{\pi_{\Iscr}(y)}$
then for any $\Jscr\st\Iscr$ one has
$\ov\Oscr_{\pi_{\Jscr}(w)}= \ov\Oscr_{\pi_{\Jscr}(y)}.$
\end{itemize}
\end{defi}
Note that by theorem \ref{2.4} the chain order on orbital varieties (respectively on $W$)
is an extension of the geometric order (that is $\Vscr\go \Wscr$ implies $\Vscr\cho \Wscr$).

\subsection{}\label{3.11}
Applying \ref{3.10} to $\pi_{i,j}$ in the case of $\gs\gl_n$ we get chain order on $\bS_n$
and Young tableaux:
\begin{defi} For $y,w \in \bS_n$ (respectively for $T,S\in \bT_n$) put $y\cho w$ (resp. $T\cho S$)
if
\begin{itemize}
\item[(i)] for any $i,j:\ 1\leq i<j\leq n$ one has $\sh(P(\pi_{i,j}(y)))\leq \sh(P(\pi_{i,j}(w)))$
(resp. $\sh(\pi_{i,j}(T))\leq \sh(\pi_{i,j}(S))$)
\item[(ii)] if for some $i,j:\ 1\leq i<j\leq n$ one has $\sh(P(\pi_{i,j}(y)))= \sh(P(\pi_{i,j}(w)))$
(resp. $\sh(\pi_{i,j}(T))= \sh(\pi_{i,j}(S))$) then for any $k,l:\ i\leq k<l\leq j$
one has $\sh(P(\pi_{k,l}(y)))= \sh(P(\pi_{k,l}(w)))$
(resp. $\sh(\pi_{k,l}(T))= \sh(\pi_{k,l}(S))$).
\end{itemize}
\end{defi}
Note that refinement (ii) is absolutely natural. We need it to sort out cases  where
two different orbital varieties associated to the same orbit are in order. In $\gs\gl_n$
such example occurs for the first time for $n=6.$ Indeed, without (ii) we will get that
$T<S$ in chain order where $T$ is from \ref{3.9}, that is
$$
T=
\vcenter{
\halign{& \hfill#\hfill
\tabskip4pt\cr
\multispan{7}{\hrulefill}\cr
\ssa
\vb& 1 & &3  & & 4 &\ts\vb\cr
\vsa
&&&&\multispan{3}{\hrulefill}\cr
\ssa
\vb& 2 & &6 & \ts\vb\cr
\vsa
&&\multispan{3}{\hrulefill}\cr
\ssa
\vb&5 & \ts\vb\cr
\vsa
\multispan{3}{\hrulefill}\cr}}\quad \longleftrightarrow\quad
\begin{array}{llllll}
(3,2,1)&(3,1,1)&(3,1)&(2,1)&(1,1)&(1)\\
(3,2)&(3,1)&(3)&(2)&(1)&\\
(2,2)&(2,1)&(2)&(1)&&\\
(2,1)&(1,1)&(1)&&&\\
(2)&(1)&&&&\\
(1)&&&&&\\
\end{array}
$$
and
$$
S=
\vcenter{
\halign{& \hfill#\hfill
\tabskip4pt\cr
\multispan{7}{\hrulefill}\cr
\ssa
\vb& 1 & &3  & & 6 &\ts\vb\cr
\vsa
&&&&\multispan{3}{\hrulefill}\cr
\ssa
\vb& 2 & &4 & \ts\vb\cr
\vsa
&&\multispan{3}{\hrulefill}\cr
\ssa
\vb&5 & \ts\vb\cr
\vsa
\multispan{3}{\hrulefill}\cr}}\quad \longleftrightarrow\quad
\begin{array}{llllll}
(3,2,1)&(2,2,1)&(2,2)&(2,1)&(1,1)&(1)\\
(3,1,1)&(2,1,1)&(2,1)&(2)&(1)&\\
(2,1,1)&(1,1,1)&(1,1)&(1)&&\\
(2,1)&(1,1)&(1)&&&\\
(2)&(1)&&&&\\
(1)&&&&&\\
\end{array}
$$

\subsection{}\label{3.12}
Note that as a straightforward corollary of definition and Gerstenhaber construction
we get that the chain order has the same 2 properties listed in proposition \ref{2.6}
as the induced Duflo order, namely
\begin{cor}
For any Young tableaux $S,T\in \bT_n$ one has
\begin{itemize}
\item[(i)] If $S\cho T$ then for any $i,j\ : 1\leq i<j\leq n$ $\pi_{i,j}(S)\cho\pi_{i,j}(T);$
\item[(ii)] $S\cho T$ iff $S^{\dagger}\chg T^{\dagger}.$
\end{itemize}
\end{cor}

Again as in \ref{2.6} let us note that $w\cho y$ iff $w_ow\chg w_oy.$

\subsection{}\label{3.13}
However, the chain order is not preserved under RS insertions.
The first examples occur in $\gs\gl_{\sr 7}.$ For example take
$$
T=
\vcenter{
\halign{& \hfill#\hfill
\tabskip4pt\cr
\multispan{7}{\hrulefill}\cr
\ssa
\vb& 1 & &2  & & 6 &\ts\vb\cr
\vsa
&&&&\multispan{3}{\hrulefill}\cr
\ssa
\vb& 3 & &5 & \ts\vb\cr
\vsa
&&&&\cr
\ssa
\vb&4 & & 7& \ts\vb\cr
\vsa
\multispan{5}{\hrulefill}\cr}}\qquad {\rm and}\qquad
S=
\vcenter{
\halign{& \hfill#\hfill
\tabskip4pt\cr
\multispan{7}{\hrulefill}\cr
\ssa
\vb& 1 & &2  & & 6 &\ts\vb\cr
\vsa
&&&&\multispan{3}{\hrulefill}\cr
\ssa
\vb& 3 & &7 & \ts\vb\cr
\vsa
&&\multispan{3}{\hrulefill}\cr
\ssa
\vb&4 &  \ts\vb\cr
\vsa
&&\cr
\ssa
\vb&5 &  \ts\vb\cr
\vsa
\multispan{3}{\hrulefill}\cr}}
$$
One can see at once that $T\chos S$. Moreover one can see at once
that $\pi_{1,6}(T)\dos \pi_{1,6}(S)$ and $\pi_{2,7}(T)\dos
\pi_{2,7}(S).$ However, recalling from \ref{2.6} $\ov T_a$ we get
$$
(5\Ra \ov T_5)=\vcenter{ \halign{& \hfill#\hfill \tabskip4pt\cr
\multispan{7}{\hrulefill}\cr \ssa \vb& 1 & &2  & & 7 &\ts\vb\cr \vsa
&&&&\multispan{3}{\hrulefill}\cr \ssa \vb& 3 & &6 & \ts\vb\cr \vsa
&&&&\cr \ssa \vb&4 & & 8& \ts\vb\cr \vsa
&&\multispan{3}{\hrulefill}\cr \ssa \vb& 5 &\ts\vb\cr \vsa
\multispan{3}{\hrulefill}\cr}}\qquad {\rm and}\qquad (5\Ra \ov S_5)=
\vcenter{ \halign{& \hfill#\hfill \tabskip4pt\cr
\multispan{7}{\hrulefill}\cr \ssa \vb& 1 & &2  & & 7 &\ts\vb\cr \vsa
&&&&&&\cr \ssa \vb& 3 & &6 && 8& \ts\vb\cr \vsa
&&\multispan{5}{\hrulefill}\cr \ssa \vb&4 &  \ts\vb\cr \vsa &&\cr
\ssa \vb&5 &  \ts\vb\cr \vsa \multispan{3}{\hrulefill}\cr}}
$$
So that $(5\Ra\ov T_5)\not\cho (5\Ra \ov S_5).$

It was shown by Barbash and Vogan that the algebraic order is preserved
under RS insertions (cf. \ref{4.5}) and we  will show in  \ref{6.6} that the geometric order
is also preserved under RS insertions. Thus, for $n\geq 7$
the chain order is a proper extension of the geometric and the algebraic orders.

\section {Primitive ideals and associated varieties}

\subsection{}\label{4.1}
Let $\ga$ denote some subalgebra of $\gog$ and let
Let $U(\ga)$ denote is universal subalgebra.

An ideal of algebra is called primitive if it is the annihilator
of some irreducible representation of this algebra.
\par
Let $\Mscr$ be the set of irreducible representations of $U(\gog)$
and $\Mscr_0$ a subset of
irreducible representations
with trivial central character. Set
$$\Xscr_0=\{Ann(M)\ :\ M\in\Mscr_0\}$$
We want to study $\Xscr_0$ as an ordered set.

Recall from \ref{1.2} and \ref{1.5} that we fix a triangular
decomposition  $\gog=\gn\bigoplus\gh\bigoplus\gn^{-}$ and denote by
$R$ the set of non-zero roots, by  $R^+$ the set of positive roots
corresponding to $\mathfrak n,$ by $\Pi\st R^+$ the resulting set of
simple roots and by $\rho$ the half-sum of positive roots. Let
$\gb=\gh\oplus\gn$ denote a Borel subalgebra of $\gog.$

For $w\in W$ let
$$M_w=U(\gog)\bigotimes_{U(\gb)}{\Co}_{-w(\rho)-\rho}$$
denote Verma module with the highest weight $-w(\rho)-\rho$
and let $L_w$ denote its (unique) simple quotient. It is called
a simple highest weight module (with the highest weight $ -w(\rho)-\rho$).
Set $I_w={\rm Ann}(L_w)$ to be the
corresponding primitive ideal in $U(\gog)$ (more precise in $U(\gn^-)$).
\par
A theorem of Duflo \cite{D} gives the
surjection from $W$ onto $\Xscr_0$ as follows
\begin{theorem} For every $I\in\Xscr_0$ there exist $w\in W$
such that $I=I_w.$
\end{theorem}
\subsection{}\label{4.2}
The surjection $\psi\,:\, W\rar \Xscr_0$ gives a decomposition of $W$ into
the left algebraic cells $\Cscr_L$ by
$$\Cscr_L(w)=\{y\in W\, :\, I_{y^{-1}}=I_{w^{-1}}\}.$$
It is customary (in the theory of primitive ideals) to call these cells simply left cells but in our context
we prefer to omit the word ``left'' and call them ``algebraic cells''
to emphasize their algebraic nature. Given an algebraic cell $\Cscr$ and some $w\in \Cscr$ we put
$I_{\Cscr}:=I_w.$
\par
We define algebraic double cell to be the union of
left cells connected via $y^{\sr -1}$:
$$\Cscr_D(w)=\bigcup_{y\in\Cscr_L(w)}\Cscr_L(y^{\sr -1}).$$

The study of $\Xscr_0$ as an ordered set can be translated into partial ordering of $W$ and of
algebraic cells. For $w,y\in W$ we put $w\ao y$ if $I_w\st I_y$ and $w\aos y$ if $I_w\subsetneq I_y$.
Respectively if $\Cscr_1, \Cscr_2$ are algebraic cells we put $\Cscr_1\ao \Cscr_2$ if
$I_{\Cscr_1}\st I_{\Cscr_2}.$

The truth of Kazhdan-Lusztig conjecture \cite{K-L} gives us a full combinatorial description of $\ao.$
We do not use Kazhdan-Lusztig combinatorics in this paper although we used it in technical
calculations explained in \ref{1.11} and it is a basis for some properties of $\ao$ we
quote here. This full combinatorial description can be found in many places
beginning from the original paper of D. Kazhdan and G. Lusztig \cite{K-L}. We will not give it here.

However, we will formulate a few related results of A. Joseph, D. Vogan and D. Barbash essential in our
further analysis.

\subsection{}\label{4.3}
We need the notion of $\tau$-invariant.
Let $w$ be any element of $W.$ Set  $S(w):=R^+ \cap^w R^-=\{\al\in R^+\ :\ w^{\sr -1}(\al)\in R^-\}.$
Set $\tau(w)=\Pi\cap S(w)$
As it is shown in \cite{B-J} and \cite{D} for primitive ideals and as it
can be seen at once from Steinberg's construction for orbital variety closures, one has
\begin{prop} Let $w,y\in W.$
\begin{itemize}
\item[(i)] If $I_{w^{-1}}\st I_{y^{-1}}$ then $\tau(w)\st\tau(y).$
\item[(ii)] If $\ov\Vscr_w\st \ov\Vscr_y$ then $\tau(w)\supset \tau(y).$
\end{itemize}
\end{prop}
In particular $\tau$-invariant is constant on
algebraic cell and on geometric cell and  we can define
$$\tau(I_{w^{-1}}):=\tau(w)\qquad {\rm and}\qquad \tau(\Vscr_w):=\tau(w).$$

\subsection{}\label{4.4}
Let us return to the case $\gog=\gs\gl_n.$ In that case $W=\bS_n$
and $\Pi=\{\al_i\}_{i=1}^{n-1}.$ Recall Robinson-Schensted procedure
$w\rar (P(w),Q(w)).$

As we have mentioned in \ref{1.6} by \cite{J-JI} and \cite{J-JIII}
one has
\begin{theorem} For $\gog=\gs\gl_n$
one has $I_{w^{-1}}=I_{y^{-1}}$ if and only if $P(w)=P(y).$
\end{theorem}
In particular by \ref{4.3} it is obvious that we should define $\tau$ invariant
on a Young tableau.
As it is mentioned in \ref{1.10} for a standard Young tableau $T\in\bT_n$ we define
$\tau(T)=\tau(T):=\{\al_i\ :\ r_{\sr T}(i+1)>r_{\sr T}(i)\}.$
As we have shown in Part I, 2.4.14 one has $\tau(P(w))=\tau(w)$ which shows that
our definition of $\tau(T)$ is consistent with other $\tau$-invariants.

Let us note also that for $\gs\gl_n$ it is very  easy to compute $\tau(w).$
Indeed, in that case, $\alpha_i\in\tau(w)$ iff $w^{-1}(i)>w^{-1}(i+1).$

\subsection{}\label{4.5}
Let us note that algebraic order has the same properties as induced
Duflo order described in \ref{2.6}.
Recall notation from \ref{2.4}.

Let $U(\gl_{\sr \Iscr})$ be the universal enveloping algebra of $\gl_{\Iscr}.$
Set $\rho_{\sr \Iscr}=0.5\sum\limits_{\alpha\in R^+_{\Iscr}}\alpha.$ For
$w\in W_\Iscr$ let $M^{\Iscr}_w$ be Verma module over $U(\gl_{\sr \Iscr})$ of the highest weight
$-w(\rho_{\sr \Iscr})-\rho_{\sr \Iscr},$ let $L_w^{\Iscr}$ be the corresponding simple quotient and
$I_w^{\Iscr}$ the corresponding primitive ideal in $U(\gl_{\Iscr}).$

Barbash and Vogan in \cite[2.24, 3.7]{B-V} provide some elementary properties of algebraic order in
any simple Lie algebra. Those are exactly the properties of induced Duflo order we have considered
in\ref{2.6}. Again we formulate them here only
for $\gog=\gs\gl_n,$ since we do it in terms of Young tableaux.
\begin{prop}
For any Young tableaux $S,T\in \bT_n$ one has
\begin{itemize}
\item[(i)] If $S\ao T$ then for any $i,j\ : 1\leq i<j\leq n$ $\pi_{i,j}(S)\ao\pi_{i,j}(T);$
\item[(ii)] If $S\ao T$ then for any $a\in\{i\}_{i=1}^{n+1}$ one has
$(\ov S_a\Downarrow a)\ao (\ov T_a\Downarrow a)$ and $(a\Ra\ov
S_a)\ao (a\Ra \ov T_a).$
\item[(iii)] $S\ao T$ iff $S^{\dagger}\ao T^{\dagger}.$
\end{itemize}
\end{prop}

Note that for the algebraic order one also has by \cite[2.24]{B-V} $w\ao y$ iff $w_ow\ag w_oy.$
\subsection{}\label{4.6}
Now we are ready to explain in detail the connection between primitive ideals and
orbital varieties described in \ref{1.6}. We return to a general semi-simple Lie algebra $\gog.$
\par
Let $M$ be a finitely
generated module over $U(\gog).$ Let $\gr\, M$
be the associated graded module over the symmetric algebra
$S(\gog)$ with respect to a good (degree) filtration on $M.$
Let $I(\gr\, M)=\Ann_{S(\gog)}(\gr\, M).$
The associated variety of $M$ is defined to be the support of $\gr M$
in $\gog^*,$ that is the variety of zeros of of $I(\gr\, M).$
$$V(M)={\rm Supp}(\gr M).$$
Identifying $\gog^*$ with $\gog$ via the Killing form we consider $V(M)$
as a subvariety of $\gog.$ If $M$ has a trivial central character then
$V(M)$ is a subvariety of the nilpotent cone $\Nscr.$

In particular consider $U(\gog)/I_w$ as a $U(\gog)$ module.
As it is shown in \cite{Jo1}, \cite{B-B I} one has the
following
\begin{theorem}
For every $w\in W$
there exist a nilpotent orbit
$\Oscr\st\Nscr$ such that $V(U(\gog)/I_w)=\ov\Oscr.$
\end{theorem}
One can refine the picture considering the associated variety of
$L_w.$ Denote by $<$  Bruhat order on $W$ (defined,  for example in \cite{Ca})
As it is shown in \cite{B-B} and \cite{Jo1} the closures of orbital varieties
are the irreducible components of associated variety of $L_w.$ Combining the information
from \cite[\S 6]{B-B}, and \cite[\S 8-9]{Jo1} in one theorem we get
\par
\begin{theorem}
For each $w\in W$ there
exists a subset $\Gamma(w)$ of $W$ such that
$$V(L_{w^-1})=\bigcup_{y\in \Gamma (w)}\ov\Vscr_y$$
where $\Gamma(w)$ has the following properties
\begin{itemize}
\item[1.] $w\in \Gamma(w)$
\item[2.] $\Gamma(w)\st\Cscr_D(w)$
\item[3.] If $y\in \Gamma(w)$ then $y\leq w$
\item[4.] If $y\in \Gamma(w)$ then $\tau(y)\supset \tau(w)$
\item[5.] $\Gamma(w^{-1})=\Gamma(w)^{-1}$
\end{itemize}
\end{theorem}
\par
As a corollary of (4) and (5) we get
\begin{itemize}
\item[6.] If $y\in \Gamma(w)$ then $\tau(y^{-1})\supset \tau(w^{-1}).$
\end{itemize}

\subsection{}\label{4.7}
The following theorem \cite[6.3]{B-B} and \cite[6.6]{Jo1}  describes
the behaviour of associated varieties on
algebraic cells.
\begin{theorem}
Let $w,y\in W.$ Then $I_{w^{\sr -1}}\st I_{y^{\sr -1}}$ implies
$V(L_w)\supset V(L_y).$
In particular, $V(L_w)$ is constant on each algebraic cell.
\end{theorem}
\subsection{}\label{4.8}
The following computations of T. Tanisaki \cite{T} show that the
associated variety of a simple highest weight module
need not be irreducible.
\par
Consider Lie algebra of type $C_3.$
Let $\Pi=\{\al_{\sr 1},\al_{\sr 2},\al_{\sr 3}\}$ be the set of fundamental
roots, where $\al_{\sr 1}$ is the long root. Set $s_i=s_{\al_i}.$
Consider $w=s_{\sr 2}s_{\sr 1}s_{\sr 2}s_{\sr 3}s_{\sr 2}s_{\sr 1}s_{\sr 2}$
and $y= s_{\sr 2}s_{\sr 3}s_{\sr 2}.$ One has $\Oscr(w)=\Oscr(y)$ and
$\Vscr_w\ne\Vscr_y.$ On the other hand
$$V(L_w)=\ov\Vscr_w\cup\ov\Vscr_y.$$
Hence $V(L_w)$ is not irreducible. This is the only non-irreducible
associated variety of a simple highest weight module in $C_3.$
\par
Now consider Lie algebra of type $B_3$ and let
$\Pi=\{\al_{\sr 1},\al_{\sr 2},\al_{\sr 3}\}$ be the set of fundamental
roots, where $\al_{\sr 1}$ is the short root. Again set $s_i=s_{\al_i}.$
Consider $w=s_{\sr 1}s_{\sr 2}s_{\sr 3}s_{\sr 1}s_{\sr 2}s_{\sr 1}$ and
$y=s_{\sr 1}s_{\sr 3}.$ Again one has $\Oscr(w)=\Oscr(y)$ and
$\Vscr_w\ne\Vscr_y.$ And again
$$V(L_w)=\ov\Vscr_w\cup\ov\Vscr_y$$
so that $V(L_w)$ is not irreducible. This is the only example
of non-irreducible associated variety in $B_3.$
\subsection{}\label{4.9}
Consider $\gog=\gs\gl_n.$
The following proposition \cite[9.12]{Jo1} is valid only for $\gs\gl_n.$
\begin{prop}
For each $w\in W,\ \ov\Vscr_w$ is the unique component
$\Vscr$ of $V(L_w)$ such that
$\tau(\Vscr)=\tau(I_{w^{\sr -1}}).$
\end{prop}
In other words for any $y\in \Gamma(w)\setminus \{w\}$ one has $\tau(y)\supsetneq \tau(I_{w^{\sr -1}}).$

\subsection{}\label{4.10}
As it is shown in \cite{M1} we get as an immediate corollary of this proposition
\begin{theorem}
For $\gog=\gs\gl_n$ one has $V(L_w)=\ov\Vscr_w.$
\end{theorem}
\Pf
Since the proof is straightforward we quote it here for completeness.

Each double cell $\Cscr_D$ is a union of finite number of left cells,
thus, there exists $\Cscr_L\st \Cscr_D$ with maximal $\tau$-invariant. (One can take
for example a left cell corresponding to some nilradical. We have at least one such
cell for any double cell.) Then by the proposition \ref{4.9} above
one has that $\Gamma(y)=\{y\}$ for any $y\in\Cscr_L.$ Thus, by \ref{4.7}[6] we get
that $\Gamma(y^{-1})=\{y^{-1}\}$ for any $y\in \Cscr_L$ which means that
$V(L_{y^{-1}})=\ov\Vscr_{y^{-1}}.$

On the other hand
by \ref{4.2} $\Cscr_D=\bigcup_{y\in\Cscr_L}\Cscr_L(y^{\sr -1}).$ Thus, for any
$w\in\Cscr_D$ there exists $y\in\Cscr_L$ such that $w\in\Cscr_L(y^{-1}).$
One has $V(L_w)=V(L_{y^{-1}})$ by \ref{4.8} and by the previous construction
$V(L_{y^{-1}})=\ov\Vscr_{y^{-1}}.$ Finally by \ref{4.4} $\ov\Vscr_{y^{-1}}=\ov\Vscr_w$
which completes the proof.
\QED

\section{ Vogan's $\Tscr_{\alpha,\beta}$ operator}

\subsection{}\label{5.1}
Let us explain Vogan's $\Tscr_{\al,\be}$
operator for primitive ideals.
Let $\al,\be\in\Pi$ be the adjacent fundamental roots of type $A_2$
i.e. such that
$s_{\al}s_{\be}s_{\al}=s_{\be}s_{\al}s_{\be}.$ We define the domain
of $\Tscr_{\al,\be}$ to be
$$D_{\al,\be}=\{w\in W\ :\ \al\not\in\tau(w),\ \be\in\tau(w)\}$$
For $w\in D_{\al,\be}$ we set
$$\Tscr_{\al,\be}(w)=
\begin{cases}
s_\al w& {\rm if}\  \be\not\in \tau(s_\al w),\\
s_\be w& {\rm otherwise.}\\
\end{cases}
$$
Note that $\Tscr_{\al,\be}(w)\in D_{\be,\al}$ and
$\Tscr_{\be,\al}(\Tscr_{\al,\be}(w))=w.$

The result of D. Vogan \cite[3.5, 3.6]{V2} gives
\begin{theorem}
For $w,y\in D_{\al,\be}$ one has
$I_{w^{-1}}\st I_{y^{-1}}$ if and only if $I_{(\Tscr_{\al,\be}(w))^{-1}}\st I_{(\Tscr_{\al,\be}(y))^-1}.$
\end{theorem}
In other words $\Tscr_{\al,\be}:D_{\al,\be}\rar D_{\be,\al}$ is an algebraic order
isomorphism.

\subsection{}\label{5.2}
Let us return to the case $\gog=\gs\gl_n.$
As a straightforward corollary of this theorem and of theorem \ref{4.4} we get that for
$y,w \in D_{\al,\be}$ one has  $P(y)=P(w)$ iff
$P(\Tscr_{\al,\be}(y))=P(\Tscr_{\al,\be}(w)).$
Moreover, since $\tau$-invariant is constant on a
cell we can define $D_{\al,\be}$ on cells as well by
$$D_{\al,\be}=\{\Cscr :\ \al\not\in\tau(\Cscr),\ \be\in\tau(\Cscr)\}$$
Note that $\Tscr_{\al,\be}(\Cscr):=\{\Tscr_{\al,\be}(y)\ |\ y\in \Cscr\}$
is well defined for any $\Cscr\in D_{\al,\be}.$

Respectively one can define $\Tscr_{\al,\be}$ also on Young tableaux. Let us give
the combinatorial description of $\Tscr_{\al,\be}(w)$ and of $\Tscr_{\al,\be}(T).$

Since $\Tscr_{\be,\al}(\Tscr_{\al,\be}(w))=w$ it is enough to consider the case
$\al=\al_i,\ \be=\al_{i+1}.$ By \ref{4.4} $w\in D_{\al_i,\al_{i+1}}$
iff $w^{-1}(i)<w^{-1}(i+1)$ and $w^{-1}(i+1)>w^{-1}(i+2)$ that is if
$\pi_{i,i+2}(w)=[i,i+2,i+1]$ or $\pi_{i,i+2}(w)=[i+2,i,i+1].$ To get $\Tscr_{\al_i,\al_{i+1}}(w)$
we interchange two entries in $w$ as follows
$$\Tscr_{\al_i,\al_{i+1}}(w)=
\begin{cases}
[\ldots,i+1,\ldots,i+2,\ldots,i,\ldots] & {\rm if}\ w=[\ldots,i,\ldots,i+2,\ldots,i+1,\ldots],\\

[\ldots,i+1,\ldots,i,\ldots,i+2,\ldots] & {\rm if}\ w=[\ldots,i+2,\ldots,i,\ldots,i+1,\ldots].\\
\end{cases}
$$

Respectively by \ref{4.4} $T\in D_{\al_i,\al_{i+1}}$ iff $r_{\sr T}(i)\geq r_{\sr T}(i+1)$
and $r_{\sr T}(i+1)< r_{\sr T}(i+2).$ Here again we have to interchange 2 entries of $T.$
If $r_{\sr T}(i)<r_{\sr T}(i+2)$ we have to interchange $i+2$ and $i+1.$
If $r_{\sr T}(i)\geq r_{\sr T}(i+2)$
we have to change $i$ and $i+1.$ Let us illustrate this by a simple example:
$$
\left\{ T= \vcenter{\halign{& \hfill#\hfill \tabskip4pt\cr
\multispan{5}{\hrulefill}\cr \ssa \vb & 1 &  &  3 &\ts\vb\cr \vsa
&&&&\cr \ssa \vb & 2 &  &  4 &\ts\vb\cr \vsa
\multispan{5}{\hrulefill}\cr}};\quad S=\vcenter{\halign{&
\hfill#\hfill \tabskip4pt\cr \multispan{5}{\hrulefill}\cr \ssa \vb &
1 &  &  3 &\ts\vb\cr \vsa &&\multispan{3}{\hrulefill}\cr \ssa \vb &
2 &\ts\vb\cr \vsa &&\cr \ssa \vb & 4 &\ts\vb\cr \vsa
\multispan{3}{\hrulefill}\cr}} \right\}\in D_{\al_2,\al_3}
$$
One has
$$
\Tscr_{\al_2,\al_3}(T)= \vcenter{\halign{& \hfill#\hfill
\tabskip4pt\cr \multispan{5}{\hrulefill}\cr \ssa \vb & 1 &  &  2
&\ts\vb\cr \vsa &&&&\cr \ssa \vb & 3 &  &  4 &\ts\vb\cr \vsa
\multispan{5}{\hrulefill}\cr}};\qquad \Tscr_{\al_2,\al_3}(S)=
\vcenter{\halign{& \hfill#\hfill \tabskip4pt\cr
\multispan{5}{\hrulefill}\cr \ssa \vb & 1 &  &  4 &\ts\vb\cr \vsa
&&\multispan{3}{\hrulefill}\cr \ssa \vb & 2 &\ts\vb\cr \vsa &&\cr
\ssa \vb & 3 &\ts\vb\cr \vsa \multispan{3}{\hrulefill}\cr}}
$$

\subsection{}\label{5.3}
Since as it is shown in \ref{4.3} $\tau$-invariant  is constant on
orbital variety we can
define $\Vscr\in D_{\al,\be}$ if $\al\not\in
\tau(\Vscr),\ \be\in\tau(\Vscr).$ Taking some $w\in W$ such that $\Vscr=\Vscr_w$
we define $\Tscr_{\al,\be}(\Vscr_w)=\Vscr_{\Tscr_{\al,\be}(w)}.$

Recall notion $\gm_{\al}$ from \ref{2.4}.
Using Vogan's calculus for orbital varieties in
$\gs\gl_n$ A. Joseph has shown in \cite[9.11]{Jo1} the following
\begin{prop}
Let $\gog=\gs\gl_n.$
If $\Vscr_w\in D_{\al,\be}$
then $\Vscr_{\Tscr_{\al,\be}(w)}$ is the unique component of
$\gm_{\al}\cap P_{\al}(\Vscr_w)$ lying in $D(\Tscr_{\be,\al}).$
\end{prop}

\subsection{}\label{5.4}
We need a very easy corollary of proposition \ref{5.3}
\begin{cor}
Let $\gog=\gs\gl_n.$
If $\Vscr_w\in D_{\al,\be}$
then $\Vscr_{\Tscr_{\al,\be}(w)}$ is the unique component of
$\gm_{\al}\cap P_{\al}(\ov{\Vscr_w})$ lying in $D_{\be,\al}.$
\end{cor}
\Pf Set $d:=\dim\Vscr_w.$ Consider $\gm_{\al}\cap\ov\Vscr_w.$ Let
$\Wscr$ be a component of this intersection. Then $\bP_{\al}\Wscr$
is a component of $\gm_{\al}\cap \bP_{\al}(\ov\Vscr_w).$ Note that
$\ov{\gm_{\al}\cap \Vscr_w}\st \gm_{\al}\cap\ov\Vscr_w$ and is
closed in it. Since $\gm_{\al}$ is a hyperplane of $\gn$ both
intersections  are equidimentional of co-dimension $1$ in
$\ov\Vscr_w.$ Let $\Wscr$ be a component of
$\gm_{\al}\cap\ov\Vscr_w$ but not a component of $\ov{\gm_{\al}\cap
\Vscr_w}.$ Then $\Wscr$ is a component of
$\gm_{\al}\cap\ov\Vscr_w\setminus\Vscr_w.$ Consider $\Oscr$ a
nilpotent orbit such that $\bG\Wscr=\ov\Oscr.$ Since
$\Wscr\st\gm_{\al}$ one has $\bP_{\al}\Wscr\st\ov\Oscr\cap\gn.$ On
the other hand $\bG\Vscr_w\cap \Wscr=\emptyset$ so $\dim \Oscr\leq
2d-2.$ Hence $\dim \bP_{\al}\Wscr\leq d-1$ and
$\bP_{\al}\Wscr=\Wscr.$ In particular $\bP_{\al}\Wscr\st\gm_{\be}$
so that $\Wscr\not\in D_{\be,\al}.$ We conclude that if
$\bP_{\al}\Wscr$ is a component of $\gm_{\al}\cap
\bP_{\al}(\ov\Vscr_w)$ lying in $D_{\be,\al}$ it must be a component
of $\ov{\gm_{\al}\cap \Vscr_w}.$ By proposition \ref{5.3}
$\ov\Vscr_{\Tscr_{\al,\be}(w)}$ is the only such component of
$\ov{\gm_{\al}\cap\Vscr_w}.$ \QED

\subsection{}\label{5.5}
Now we are ready to show that $\Tscr_{\al,\be}$ is a geometric order isomorphism.
\begin{theorem}
Let $\gog=\gs\gl_n.$
For $w,y\in D_{\al,\be}$ one has
$\Vscr_w\st \ov\Vscr_y$ if and only if $\Vscr_{(\Tscr_{\al,\be}(w))}\st \Vscr_{(\Tscr_{\al,\be}(y))}.$
\end{theorem}
\Pf
By symmetry of $\Tscr_{\al,\be}$ it is enough to show
only one direction. $\Vscr_w\st\ov\Vscr_y$ implies that
$\ov{\gm_{\al}\cap \bP_{\al}(\Vscr_w)}\st
\gm_{\al}\cap \bP_{\al}( \ov\Vscr_y).$
By  proposition \ref{5.3}
$\Vscr_{\Tscr_{\al,\be}(w)}$ is a component of
$\ov{\gm_{\al}\cap \bP_{\al}(\Vscr_w)}$ and by corollary \ref{5.4}
$\ov\Vscr_{\Tscr_{\al,\be}(y)}$ is the only component of
$\gm_{\al}\cap \bP_{\al}(\ov\Vscr_y)$ not lying in $\gm_{\be}.$
Hence by irreducibility of $\Vscr_{\Tscr_{\al,\be}(w)}$
it must lie in $\ov\Vscr_{\Tscr_{\al,\be}(y)}.$
\QED

\subsection{}\label{5.6}
By \ref{5.1} and \ref{5.5} $\Tscr_{\al,\be}$ preserves both algebraic and geometric
order. Straightforward checking shows that $\Tscr_{\al,\be}$ preserves induced Duflo
order for $n\leq 5.$ Moreover the induced Duflo order
coincides with the chain order for $n\leq 5,$ therefore coincides with
the algebraic and geometric orders which are sandwiched between the induced
Duflo and the chain orders.

However, for $n\geq 6$\ $\Tscr_{\al,\be}$ does not preserve
induced Duflo order anymore, so that induced Duflo order is a proper restriction
the algebraic order for $n\geq 6.$

Let us show this by the example.
Consider $w=[3,\, 5,\, 6,\,1,\,2,\,4]$ and $ws_{\sr 1}=[5,\, 3,\, 6,\,1,\,2,\,4].$
One has that $w\dor ws_{\sr 1}.$ Let $T:=P(w)$ and $S:=P(ws_{\sr 1}).$
One has $T\dor S$ where by Robinson-Schensted procedure.
$$T=\vcenter{
\halign{&\hfill#\hfill
\tabskip4pt\cr
\multispan{7}{\hrulefill}\cr
\ssa
\vb &1 & & 2  & & 4 & \ts\vb\cr
\vsa
&&&&&&\multispan{3}{\hrulefill}\cr
\ssa
\vb & 3 && 5 && 6  &  \ts\vb\cr
\vsa
\multispan{7}{\hrulefill}\cr}},\qquad
S=\vcenter{
\halign{&\hfill#\hfill
\tabskip4pt\cr
\multispan{7}{\hrulefill}\cr
\ssa
\vb &1 & & 2  && 4 & \ts\vb\cr
\vsa
&&&&\multispan{3}{\hrulefill}\cr
\ssa
\vb & 3 && 6 &  \ts\vb\cr
\vsa
&&\multispan{3}{\hrulefill}\cr
\ssa
\vb & 5&  \ts\vb\cr
\vsa
\multispan{3}{\hrulefill}\cr}}.$$

By \ref{4.5} $\tau(T)=\tau(S)=\{\al_{\sr 2},\, \al_{\sr 4}\},$ thus, $T,S\in D_{\al_3,\al_4}.$
By \ref{5.2}
$$\Tscr_{\al_{\sr 3}\al_{\sr 4}}(T)=\vcenter{
\halign{&\hfill#\hfill
\tabskip4pt\cr
\multispan{7}{\hrulefill}\cr
\ssa
\vb &1 & & 2  && 3 & \ts\vb\cr
\vsa
&&&&&&\multispan{3}{\hrulefill}\cr
\ssa
\vb & 4 && 5 && 6  &  \ts\vb\cr
\vsa
\multispan{7}{\hrulefill}\cr}},\qquad
\Tscr_{\al_{\sr 3}\al_{\sr 4}}(S)=\vcenter{
\halign{&\hfill#\hfill
\tabskip4pt\cr
\multispan{7}{\hrulefill}\cr
\ssa
\vb &1 & & 2  && 5 & \ts\vb\cr
\vsa
&&&&\multispan{3}{\hrulefill}\cr
\ssa
\vb & 3 && 6 &  \ts\vb\cr
\vsa
&&\multispan{3}{\hrulefill}\cr
\ssa
\vb & 4&  \ts\vb\cr
\vsa
\multispan{3}{\hrulefill}\cr}}.$$
By theorem \ref{5.1} $\Tscr_{\al_3,\al_4}(T)\ao\Tscr_{\al_3,\al_4}(S)$
(thus, also  $\Tscr_{\al_3,\al_4}(T)\go\Tscr_{\al_3,\al_4}(S)$). On the other hand
as it is shown in in Part I, 4.1.6 $\dot T\not\dor \dot S.$ Note that
$\Tscr_{\al_{\sr 3}\al_{\sr 4}}(T)$ is the translation of $\dot T$ into the standard
form and $\Tscr_{\al_{\sr 3}\al_{\sr 4}}(S)$ is the translation of $\dot S$
into the standard form, thus, $\Tscr_{\al_{\sr 3}\al_{\sr 4}}(T)\not\dor \Tscr_{\al_{\sr 3}\al_{\sr 4}}(S).$

\subsection{}\label{5.7}
As well straightforward checking shows that $\Tscr_{\al,\be}$ preserves the chain order for
$n\leq 6$ and moreover the chain order coincides with the algebraic (hence also with
the geometric) order.

Let us show that  $\Tscr_{\al,\be}$ does not preserve chain order in
$\gs\gl_n$ where $n\geq 7$ so that for $n\geq 7$ the chain order is
a proper extension of the geometric order.

Indeed, consider again the example from \ref{3.13}
$$
T=
\vcenter{
\halign{& \hfill#\hfill
\tabskip4pt\cr
\multispan{7}{\hrulefill}\cr
\ssa
\vb& 1 & &2  & & 6 &\ts\vb\cr
\vsa
&&&&\multispan{3}{\hrulefill}\cr
\ssa
\vb& 3 & &5 & \ts\vb\cr
\vsa
&&&&\cr
\ssa
\vb&4 & & 7& \ts\vb\cr
\vsa
\multispan{5}{\hrulefill}\cr}}\qquad {\rm and}\qquad
S=
\vcenter{
\halign{& \hfill#\hfill
\tabskip4pt\cr
\multispan{7}{\hrulefill}\cr
\ssa
\vb& 1 & &2  & & 6 &\ts\vb\cr
\vsa
&&&&\multispan{3}{\hrulefill}\cr
\ssa
\vb& 3 & &7 & \ts\vb\cr
\vsa
&&\multispan{3}{\hrulefill}\cr
\ssa
\vb&4 &  \ts\vb\cr
\vsa
&&\cr
\ssa
\vb&5 &  \ts\vb\cr
\vsa
\multispan{3}{\hrulefill}\cr}}
$$
As we have shown $T\cho S.$ Note that $T,S\in D_{\al_5,\al_6}$ so that
$$\Tscr_{\al_5,\al_6}(T)= \vcenter{
\halign{& \hfill#\hfill
\tabskip4pt\cr
\multispan{7}{\hrulefill}\cr
\ssa
\vb& 1 & &2  & & 7 &\ts\vb\cr
\vsa
&&&&\multispan{3}{\hrulefill}\cr
\ssa
\vb& 3 & &5 & \ts\vb\cr
\vsa
&&&&\cr
\ssa
\vb&4 & & 6& \ts\vb\cr
\vsa
\multispan{5}{\hrulefill}\cr}}\qquad {\rm and}\qquad
\Tscr_{\al_5,\al_6}(S)=
\vcenter{
\halign{& \hfill#\hfill
\tabskip4pt\cr
\multispan{7}{\hrulefill}\cr
\ssa
\vb& 1 & &2  & & 5 &\ts\vb\cr
\vsa
&&&&\multispan{3}{\hrulefill}\cr
\ssa
\vb& 3 & &7 & \ts\vb\cr
\vsa
&&\multispan{3}{\hrulefill}\cr
\ssa
\vb&4 &  \ts\vb\cr
\vsa
&&\cr
\ssa
\vb&6 &  \ts\vb\cr
\vsa
\multispan{3}{\hrulefill}\cr}}
$$
One can see at once that $\sh(\pi_{1,6}(\Tscr_{\al_5,\al_6}(T)))\not<
\sh(\pi_{1,6}(\Tscr_{\al_5,\al_6}(S))).$

\section{\bf Inclusion of orbital varieties and Duflo-Vogan and chain-Vogan orders}

\subsection{}\label{6.1}
In this section we concentrate on the study of orbital variety closures.

We begin with the consideration of 3 properties of algebraic order described in \ref{4.5}, which are
also true for induced Duflo orders it is shown in \ref{2.6}.

The first property, namely, that for any reductive $\gog$ the inclusion of orbital variety closures
implies the inclusion of their projections on Levi factors, was shown to be true in Part I, \S 4.1.1.

The third property, namely, $T<S$ iff $S^{\dagger}< T^{\dagger},$ is very natural for any
combinatorially defined order (in particular for
induced Duflo and chain orders). For algebraic order it is easily shown to be true with the help of
Kazhdan-Lusztig data which is purely combinatoric.
I am sure that this also should be true for geometric order,
however, it demands more advanced combinatorial tools then those we have
at hand now.

Now we are going to show that geometric order has also  the second property,
namely, we will construct the embeddings
from Levi factors to $\gog$ preserving the inclusion of orbital variety closures.

Recall the notation from \ref{2.4}. Given $\Iscr\st \Pi$ recall that $\Vscr^{\Iscr}$ denote an orbital
variety in $\gl_{\sr \Iscr}.$ For any $f\in F_{\Iscr}$ we can define embedding
$\eps_f :\gn_{\sr \Iscr}\ha \gn$
via $\eps_f(X_{\al})=X_{f(\al)}$ for $\al\in R_{\Iscr}^+.$ Our aim is to show that these embeddings
preserve inclusions of orbital variety closures. This can be formulated as follows
\begin{theorem}
Let $\gog$ be some semi-simple Lie algebra and let $\Iscr$ be some subset of $\Pi.$
For any $y,w \in W_\Iscr$ such that $\ov\Vscr^{\Iscr}_y\st \ov\Vscr^{\Iscr}_w$ and for any
$f\in F_{\Iscr}$ one has $\ov\Vscr_{fy}\st \ov\Vscr_{fw}.$
\end{theorem}
We prove the theorem via 2 technical lemmas below.

\subsection{}\label{6.2}
 We begin with
\begin{lemma}
For any $w\in W_\Iscr$ and $f\in F_\Iscr$ one has
$$\gn\cap^{fw}\gn=f(\gn_{\sr \Iscr}\cap^w\gn_{\sr \Iscr})\oplus \gn\cap^f\gm_{\sr \Iscr}.$$
\end{lemma}
\Pf
Let us first show that $f(\gn_{\sr \Iscr}\cap^w\gn_{\sr \Iscr})\st \gn\cap^{fw}\gn,$\
$\gn\cap^f\gm_{\sr \Iscr}\st \gn\cap^{fw}\gn$ and their sum is a direct sum.

To show that $\gn\cap^f\gm_{\sr \Iscr}\st \gn\cap^{fw}\gn$
it is enough to show that for any $X_{\al}\st \gn\cap^f\gm_{\sr \Iscr}$ one has
$(fw)^{\sr -1}(X_\al)\st \gn.$ Indeed, $X_{\al}\st \gn\cap^f\gm_{\sr \Iscr}$
means that $X_\al\st \gn$ and $f^{\sr -1}(X_\al)=X_\be\st\gm_{\sr \Iscr}.$ Then
$w^{\sr -1}f^{\sr -1}(X_\al)=w^{\sr -1}(X_\be)\st \gm_{\sr \Iscr}$ since $\gm_{\sr \Iscr}$
is stable under the action of $W_\Iscr.$

Now let us show that $f(\gn_{\sr \Iscr}\cap^w\gn_{\sr \Iscr})\st \gn\cap^{fw}\gn.$
Indeed, $f(\gn_{\sr \Iscr})\st\gn$ by the definition of $F_\Iscr,$ thus,
$f(\gn_{\sr \Iscr}\cap^w\gn_{\sr \Iscr})\st \gn.$ As well
$f(\gn_{\sr \Iscr}\cap^w\gn_{\sr \Iscr})\st fw(\gn).$ Therefore
$f(\gn_{\sr \Iscr}\cap^w\gn_{\sr \Iscr})\st \gn\cap^{fw}\gn.$

Further, since $\gn_{\sr \Iscr}\cap\gm_{\sr\Iscr}=\emptyset$ we get that the sum
is direct.

It remains to show that $\gn\cap^{fw}\gn \st
f(\gn_{\sr \Iscr}\cap^w\gn_{\sr \Iscr})\oplus \gn\cap^f\gm_{\sr \Iscr}.$

Assume that $X_{\al}\st\gn\setminus {f(\gn_{\sr \Iscr}\cap^w\gn_{\sr \Iscr})\oplus \gn\cap^f\gm_{\sr \Iscr}}.$ We must show that $(fw)^{\sr -1}(X_\al)\st \gn^-.$
Indeed, if $X_{\al}\st f(\gn_{\sr \Iscr})$ then $X_\al=f(X_\be)$ where $w^{\sr -1}(X_\be)\st\gn^-$ so that
$w^{\sr -1}f^{\sr -1}(X_\al)=w^{\sr -1}(X_\be)\st \gn^-.$ Otherwise, $X_\al \st \gn\setminus
f(\gn_{\Iscr}\oplus\gm_{\Iscr})$
thus, $f^{\sr -1}(X_\al)\not\st \gn_{\sr \Iscr}$ and $f^{\sr -1}\not\st \gm_{\sr\Iscr},$ as well. Hence,
$f^{\sr -1}(X_\al)\st \gn^-$ so that again $w^{\sr -1}f^{\sr -1}(X_\al)\st \gn^-.$
\QED

\subsection{}\label{6.3}
Recall notion of $\bB_\Iscr$ from \ref{2.4}.
Exactly in the same way as in \ref{6.1} we have defined
$\eps_f :\gn_{\sr \Iscr}\ha \gn$ for some $\Iscr\st\Pi$ and any $f\in F_\Iscr$ we define
the embedding $\eps_f:\bB_\Iscr\ha \bB$ via $\eps_f(\bB_\al)=\bB_{f(\al)}$ for any $\al\in R_\Iscr^+.$

For any $w\in W$ we fix a representative in $\bG$ and denote it by $g_w.$ Recall that for any
$B\in \bB$ one has $w(B)=g_wBg_w^{-1}$ which gives us
\begin{lemma} Let $\Iscr\st \Pi$ and let $f$ be some element of $F_\Iscr.$ Then
\begin{itemize}
\item[(i)] $\eps_f:\bB_\Iscr\ha \bB$ is a monomorphism.
\item[(ii)] $\eps(\bB_\Iscr)(\gn\cap^f\gm_{\sr\Iscr})=\gn\cap^f\gm_{\sr\Iscr}.$
\end{itemize}
\end{lemma}
\Pf
(i)Indeed, this is a homomorphism since for any $A,B\in \bB_\Iscr$ one has
$$\eps_f(AB^{-1})=g_f(AB^{-1})g_f^{-1}=\eps_f(A)\eps_f(B^{-1}).$$
As well $\eps_f$ is injective by definition.

(ii) Obviously $\gn\cap^f\gm_{\sr\Iscr}\st \eps(\bB_\Iscr)(\gn\cap^f\gm_{\sr\Iscr}).$
The other inclusion is almost straightforward. Let us check it.
 Every $X\in \gn\cap^f\gm_{\sr\Iscr}$ can be represented as $X=g_fYg_f^{-1}$
for some $Y\in\gm_{\sr\Iscr}$
and every
$B\in \eps(\bB_\Iscr)$ can be represented as $B=g_fAg_f^{-1}$ for some $A\in \bB_\Iscr.$
Thus, $BXB^{-1}=g_fAYA^{-1}g_f^{-1}.$ Since $\gm_{\sr \Iscr}$ is $\bB$ stable it is also
$\bB_\Iscr$ stable. Thus, $AYA^{-1}\in \gm_{\sr \Iscr}$ and $BXB^{-1}\in f(\gm_{\sr \Iscr}).$
On the other hand $BXB^{-1}\in\gn.$ Therefore, $BXB^{-1}\in \gn\cap^f\gm_{\sr\Iscr}.$
This provides the other inclusion.
\QED

\subsection{}\label{6.4}
Now the proof of theorem \ref{6.1} follows straightforwardly.

\Pf
Assume that $w,y\in W_\Iscr$ are such that $\Vscr^\Iscr_y\st\ov{\Vscr^\Iscr_w}.$
This means that
$\gn_{\sr \Iscr}\cap^y \gn_{\sr \Iscr}\st \ov{\bB_\Iscr(\gn_{\sr \Iscr}\cap^y \gn_{\sr \Iscr})}.$
To prove the theorem it  suffices to show that for any $f\in F_\Iscr$ one has
$\gn\cap{fy}\gn\st \ov{\bB(\gn\cap{fw}\gn)}.$
Indeed,
$$\begin{array}{rll}
\gn\cap^{fy}\gn&=\eps_f(\gn_{\sr \Iscr}\cap^y \gn_{\sr \Iscr})\oplus\gn\cap^f\gm_{\sr\Iscr}&
{\rm by\ lemma\ \ref{6.2},}\\
&\st \eps_f(\ov{\bB_\Iscr(\gn_{\sr \Iscr}\cap^w \gn_{\sr \Iscr})})+\gn\cap^f\gm_{\sr\Iscr}&{\rm by\ hypothesis,}\\
&\st \ov{\eps_f(\bB_\Iscr(\gn_{\sr \Iscr}\cap^w \gn_{\sr \Iscr}))}+\gn\cap^f\gm_{\sr\Iscr}&{\rm by\
continuity \ of\ } \eps_f,\\
&\st \ov{\eps_f(\bB_\Iscr)(\eps_f(\gn_{\sr \Iscr}\cap^w \gn_{\sr \Iscr})\oplus\gn\cap^f\gm_{\sr\Iscr})}
             &{\rm by\ lemma\ }\ref{6.3},\\
&=\ov{\eps_f(\bB_\Iscr)(\gn\cap^{fw}\gn)}& {\rm by\ lemma\ \ref{6.2},}\\
&\st\ov{\bB(\gn\cap^{fw}\gn)}. & \\
\end{array}$$
\QED

\subsection{}\label{6.5}
Recall notion of $D_{\al,\be}$ from \ref{5.1}.
As an immediate corollary of theorem \ref{6.1} we get that if $w^{-1}\st D_{\al,\be}$
then $\Vscr_{(\Tscr_{\al,\be}(w^{-1}))^{-1}}=\Vscr_w.$ Indeed, consider $\Iscr={\al,\be}.$
One has $w=f_{\sr\Iscr}w\pr$ and
$(\Tscr_{\al,\be}(w^{-1}))^{-1}=f_{\sr\Iscr}(\Tscr_{\al,\be}((w\pr)^{-1}))^{-1}.$
Since $\Vscr^\Iscr_{w\pr}=\Vscr^{\Iscr}_{(\Tscr_{\al,\be}((w\pr)^{-1}))^{-1}}$ the result
follows.

This result is well know. Moreover the cells in $W$ of $\gs\gl_n$ are built exactly by the
union of these relations.

Moreover the algebraic cells in $B_n$ and $C_n$ are exactly defined by decomposition of
$B_2$ and $A_2$ using this procedure as it was shown by \cite{G}.

We can try to make the same for geometric cells in $W$ of $\gog$ of type $B_n$ or $C_n.$
We compute geometric cells of type $A_2$ and $B_2$. Then applying theorem \ref{4.1} to all
possible $\Iscr=\{\al,\be\}$ of type $A_2$ and $B_2$ we get the decomposition of $W$ into the subsets
such that each geometric cell is a union of these subsets. Unfortunately, in case of $B_n$
and $C_n$ there are geometric cells which are union of a few such subsets.

Consider for example $\gog$ of type $B_3.$ Let $\al_1$ be the short root and $\al_2,\al_3$
be the long roots. Set $s_i=s_{\al_i.}$ As it is shown in \cite{T} $\Cscr_{s_1s_2s_1}$ is
a union of 3 such subsets:
$$\Cscr_{s_1s_2s_1}=\{s_1s_2s_1\}\coprod\{s_1s_2s_1s_3,\ s_1s_2s_1s_3s_2\}\coprod\{s_1s_2s_1s_3s_2s_1\}$$

\subsection{}\label{6.6}
For $\gog=\gs\gl_n$  theorem \ref{6.1} provides us that
geometric order like algebraic and induced Duflo orders is preserved under insertions, namely
\begin{prop}
For any Young tableaux $S,T\in \bT_n$ if $S\go T$ then for any
$a\in\{i\}_{i=1}^{n+1}$ one has $(\ov S_a\Downarrow a)\go (\ov
T_a\Downarrow a)$ and $(a\Ra\ov S_a)\go (a\Ra \ov T_a).$
\end{prop}
\Pf
Since the proofs of the first and the second implications are exactly the same we show only
the first one.

Let $S,T\in \bT_n$ be such that $S\go T$ Let $w,y \in \bS_n$ be such that
$S=P(w)$ and $T=P(y).$ This means $\ov\Vscr_y\st\ov\Vscr_w.$

Recall the notion of $[\ov w_a]$ from \ref{2.6}. One has $[\ov w_a,
a]=f_aw$ and $[\ov y_a,a]=f_ay$ where $f_a=s_as_{a+1}\ldots s_n$ is
in $F_{\{a_i\}_{i=1}^{n-1}}.$ Therefore by theorem \ref{6.1} $(\ov
S_a \Da a)=P([\ov w_a, a])\go P([\ov y_a,a])=(\ov T_a\Da a).$ \QED

\subsection{}\label{6.7}
As we have shown in \S 5 both geometric and algebraic orders are preserved under
$\Tscr_{\al,\be}$ procedure and both induced Duflo order and geometric order are not
preserved under this procedure. The natural idea is to use this procedure  to strengthen
induced Duflo order on one hand and
to refine chain on the other hand.

Let us call an order generated by procedures $\Tscr_{\al,\be}(T)$ and Robinson-Schensted
insertions $(a\Ra T),\ (T\Da a)$ Duflo-Vogan order on $\bT_n$ and denote it by $\dvo.$

Explicitly for $T,S\in \bT_n$ put $T\dvo S$ if there exists a chain
of $T_1,\ldots, T_k=T$ and $S_1,\ldots,S_k=S$ such that $T_1=(\ov
T_a\pr\Da a)$ and $S_1=(\ov S_a\pr\Da a)$ or $T_1=(a\Ra \ov T_a\pr)$
and $S_1=(a\Ra \ov S_a\pr)$ where $T\pr,S\pr\in \bT_{n-1}$ are such
that $T\pr\dvo S\pr$ and for any $i\ :\ 2\leq i\leq k$ there exists
$\Tscr_{\al_{i-1},\be_{i-1}}$ such that
$T_i=\Tscr_{\al_{i-1},\be_{i-1}}(T_{i-1})$ and
$S_i=\Tscr_{\al_{i-1},\be_{i-1}}(S_{i-1}).$

Obviously one has the following: $\dor\prec\dvo\preceq \ao.$

Let us call an order generated by ordering of Young diagrams for chains and restricted
by $\Tscr_{\al,\be}$ Chain-Vogan order on $\bT_n$ and denote it by $\chvo.$

Explicitly, for $T,S\in \bT_n$ put $T\chvos S$ if $\sh(T)<\sh(S),$
$\pi_{1,n-1}(T)\chvo \pi_{1,n-1}(S),$ $\pi_{2,n}(T)\chvo \pi_{1,n-1}(S)$ and for any pair of tableaux
in a chain
$T=T_1,\ldots, T_k$ and $S=S_1,\ldots,S_k$ where for any $i\ :\ 2\leq i\leq k$
there exists $\Tscr_{\al_{i-1},\be_{i-1}}$ such that
$T_i=\Tscr_{\al_{i-1},\be_{i-1}}(T_{i-1})$ and $S_i=\Tscr_{\al_{i-1},\be_{i-1}}(S_{i-1})$
one has  $\pi_{1,n-1}(T_i)\chvo \pi_{1,n-1}(S_i)$ and $\pi_{2,n}(T_i)\chvo \pi_{1,n-1}(S_i).$

Again one has $\go\preceq \chvo\prec\cho.$

Since both Duflo-Vogan and chain-Vogan orders are of combinatorial nature, the program on Mathematica
was written to compare them. With the help of this program I have found that these orders coincide for
$n\leq 9.$

The first pair of tableaux such that $T\dvonot S$ however, $T\chvos S$ occurs
in $\bT_{10}.$ They are
$$
T=
\vcenter{
\halign{& \hfill#\hfill
\tabskip4pt\cr
\multispan{9}{\hrulefill}\cr
\ssa
\vb& 1 & &2  & & 5 && 6&\ts\vb\cr
\vsa
&&&&&&\multispan{3}{\hrulefill}\cr
\ssa
\vb& 3 & &4 &&9& \ts\vb\cr
\vsa
&&&&\multispan{3}{\hrulefill}\cr
\ssa
\vb&7 &&8&  \ts\vb\cr
\vsa
&&\multispan{3}{\hrulefill}\cr
\ssa
\vb&10 &  \ts\vb\cr
\vsa
\multispan{3}{\hrulefill}\cr}}\qquad {\rm and}\qquad
S=
\vcenter{
\halign{& \hfill#\hfill
\tabskip4pt\cr
\multispan{7}{\hrulefill}\cr
\ssa
\vb& 1 & &2  & & 6&\ts\vb\cr
\vsa
&&&&&&\cr
\ssa
\vb& 3 & &4 &&8& \ts\vb\cr
\vsa
&&&&\multispan{3}{\hrulefill}\cr
\ssa
\vb&5 && 9 &  \ts\vb\cr
\vsa
&&&&\cr
\ssa
\vb& 7 && 10 &  \ts\vb\cr
\vsa
\multispan{5}{\hrulefill}\cr}}
$$
All other examples in $\bT_{10}$ are obtained from these two tableaux by
$\Tscr_{\al,\be}$ operations and transposing $T\rightarrow T^{\dagger},$ that is
up to these operations  this is the only case in $\bT_{10}.$

To check the situation with the inclusion of orbital variety
closures $\ov\Vscr_T,\ \ov\Vscr_S$ I used the program of F. Du Cloix
for the computations of Kazhdan-Lusztig polynomials. The
computations show that corresponding primitive ideals $I_T\subset
I_S$ so that for $n\leq 10$ chain-Vogan order coincides with the
algebraic order (thus, also with the geometric order). That is,
$\ov\Vscr_S\subset\ov\Vscr_T.$

\subsection{}
\label{6.8}
As we noticed in \ref{3.13} chain order unlike  induced Duflo, Duflo-Vogan,
algebraic and geometric orders is not preserved under the insertions. As for chain-Vogan order I do not
know meanwhile whether it is preserved under the insertions for $n\geq 11$ or those insertions can be used
for the further refinement of the order. This is a very interesting question, however, we leave it for the
future research.

\section{Some properties of the cover for geometric order}
\subsection{}
\label{7.1}
Recall the notion of the cover for a partial order from \ref{1.12}.
In this last section we will describe some properties of the cover of Young tableaux
for the geometric order. We will call it in short a geometric cover of $T$ or of $\Vscr_T.$
All our results here are true also for the algebraic order. In this section ``the cover'' will mean
``the cover for the geometric order''.

We formulate everything for projection
$\pi_{{\sr 1},n-\sr 1}:\bT_n\rar \bT_{n-1}$
and induction $\Da:\bT_n\rar \bT_{n+1}$ but it can as well be formulated
for $\pi_{{\sr 2},n}$ and $\Rightarrow.$

\subsection{}
\label{7.2}
Gerstenhaber's construction provides a simple and nice description
of the cover of order on nilpotent orbits defined by the inclusion of
the closures. Let us describe it in the terms of corresponding partitions.
Let $\lambda=(\lambda_1,\ldots,\lambda_j,0)$ be some partition of $n.$
Its cover for the order defined in\ref{3.1} is constructed as follows.
\begin{itemize}
\item[(i)] For any $i:\ 1\leq i\leq j$ such that $\lambda_i\geq \lambda_{i+1}+2$
there exists  $\mu=(\mu_1,\ldots,\mu_{j+1})$ in the cover where
$\mu_s=\lambda_s$ for any $s\ne i,i+1$
and $\mu_i=\lambda_i-1,$ $\mu_{i+1}=\lambda_{i+1}+1.$
\item[(ii)] for any $i:\ 1\leq i<j$ such that $\lambda_{i+1}=\ldots=\lambda_{i+k-1}=\lambda_i-1$
and $\lambda_{i+k}=\lambda_i-2$ for some $k\geq 2$ there exists $\mu=(\mu_1,\ldots,\mu_{j+1})$
in the cover where $\mu_s=\lambda_s$ for any $s\ne i,i+k$
and $\mu_i=\mu_{i+k}=\lambda_i-1.$
\end{itemize}

Let $\Oscr\pr$ be a nilpotent orbit in the cover of nilpotent orbit
$\Oscr.$ Let $\Vscr$ be any orbital variety associated to $\Oscr.$ Then
as it is shown in Part I 4.1.8 there exists an orbital variety $\Vscr\pr$
associated to $\Oscr\pr$ such that $\Vscr\pr\gegs \Vscr.$ Obviously it is
in the cover of $\Vscr.$

However, if $\Vscr_{T\pr}\gegs\Vscr_T$ is in the cover
this does not imply that $\sh(T\pr)$ is in the cover of $\sh(T).$ The first
``jump'' (that is $S$ is the  cover of $T$, but $\sh(S)$ is not in the cover of $\sh(T)$) occurs
already in $\bT_4.$ Let us consider this example in detail.
$$
T =
\vcenter{
\halign{& \hfill#\hfill
\tabskip4pt\cr
\multispan{7}{\hrulefill}\cr
\ssa
\vb & 1 &  & 2 &  & 3 &\ts\vb\cr
\vsa
&&\multispan{5}{\hrulefill}\cr
\ssa
\vb & 4  &\ts\vb\cr
\vsa
\multispan{3}{\hrulefill}\cr}}\quad
\Cscr_T=\{s_{\sr 3},\ s_{\sr 2}s_{\sr 3},\ s_{\sr 1}s_{\sr 2} s_{\sr 3}\},$$
$$
S=
\vcenter{
\halign{& \hfill#\hfill
\tabskip4pt\cr
\multispan{5}{\hrulefill}\cr
\ssa
\vb & 1 & & 2 &\ts\vb\cr
\vsa
&&\multispan{3}{\hrulefill}\cr
\ssa
\vb & 3  &\ts\vb\cr
\vsa
\multispan{0}{\hrulefill}\cr
\ssa
\vb  & 4 &\ts\vb\cr
\vsa
\multispan{3}{\hrulefill}\cr}}\quad
\Cscr_S=\{s_{\sr 2} s_{\sr 3} s_{\sr 2}, \ s_{\sr 1} s_{\sr 2}s_{\sr 3} s_{\sr 2},
\ s_{\sr 1} s_{\sr 2}s_{\sr 1}s_{\sr 3}s_{\sr 2}\}$$
The intermediate nilpotent orbit has just two orbital varieties labeled by
$$ P=
\vcenter{
\halign{& \hfill#\hfill
\tabskip4pt\cr
\multispan{5}{\hrulefill}\cr
\ssa
\vb & 1 &  & 2 &\ts\vb\cr
\vsa
&&&&\multispan{0}{\hrulefill}\cr
\ssa
\vb & 3& &4  &\ts\vb\cr
\vsa
\multispan{5}{\hrulefill}\cr}},\qquad
Q=
\vcenter{
\halign{& \hfill#\hfill
\tabskip4pt\cr
\multispan{5}{\hrulefill}\cr
\ssa
\vb & 1 &  & 3 &\ts\vb\cr
\vsa
&&&&\multispan{0}{\hrulefill}\cr
\ssa
\vb & 2& &4  &\ts\vb\cr
\vsa
\multispan{5}{\hrulefill}\cr}}$$
These satisfy $\tau(P)\not \supset \tau(T)$
and $\tau(S)\not\supset \tau(Q).$ Hence
$T\not\gos P,\ Q\not\gos S,$ so that $S$ is the cover of $T.$

\subsection{}
\label{7.3}
The lemma below implies that this is a general phenomenon.

To formulate and prove the lemma we need to recall some combinatorial notation from
Part I, 2.4.2. Given $T\in \bT_n$ of shape $\sh(T)=(\lambda_1,\ldots\lambda_j).$
For any $p,r\ :\ 1\leq p\leq r\leq j$ put $T^{p,r}$ to be a tableau consisting of rows $p,\ldots,r$
of tableau $T.$ Note that for any $p\ : 1<p\leq j$ one has that if $x$ is a word such that
$P(x)=T^{p,j}$ and $y$ is a word such that $P(y)=T^{1,p-1}$ then as it is shown in Part I 3.2.3(v)
$P([x,y])=T.$
Recall notation $r_{\sr T}(j)$ from \ref{1.10}. Note that the information $r_{\sr T}(j)$
for all $j\ :\ 1\leq j\leq n$ determines $T$ completely, since the numbers increase in
the rows from left to right.
\begin{lemma}
Consider $T \in \bT_n$ with $\sh(T)=(\lambda_1,\ldots,\lambda_j,0)$  and
assume that $r_{\sr T}(n)=i.$ Assume that $k\geq 1$ is the minimal such that
$\lambda_{i+k}\leq \lambda_i-2.$ Then $S$ which is obtained from $T$ by moving
the box with $n$ from row $i$ to row $i+k,$ that is such that $r_{\sr S}(j)=r_{\sr T}(j)$
for any $j<n$ and $r_{\sr S}(n)=i+k$ is in the cover of $T.$
\end{lemma}
\Pf
(1)\ \  Let us show first that $T\dos S$ and, thus, also $T\gos S.$  Indeed, set
$w=[x,n,y]$ where $P(x)=T^{i+k-1,j}$ and $P(y)=(\pi_{1,n-1}(T))^{1,i+k-2}.$
Let also $w^{-1}(n)=m.$ Then by Part I, 3.2.3(v) and 3.2.5
$T=P(w),\ S=P(ws_{m-\sr 1}).$ Given $w=[a_1,\ldots,a_n],$ a well known fact that
$ws_i\dgs w$ iff $a_{i+1}>a_i$ is shown for example in \cite[pp.73-74]{Ja}. In our case this provides
$ws_{m-\sr 1} \dgs w$. Hence $T\dos S.$
\par
(2)\ \  We must show that $S$ is a geometric descendant of $T$.

We distinguish the following
two cases.
\begin{itemize}
\item[(a)] Either \ $k=1$ or $k>1$ and $\lambda_{i+k}=\lambda_i-2.$
\item[(b)] $k>1$ and $\lambda_{i+k}<\lambda_i-2.$
\end{itemize}
\par
In case (a) $S$ is a geometric  of $T$ because $\sh(S)$ is a
descendant of $\sh(T).$
\par
In case (b) $\sh(S)=(\lambda_1,\ldots,\lambda_i-1,\lambda_{i+1},\ldots,\lambda_{i+k-1},
\lambda_{i+k}+1,\ldots)$ and there exists a unique intermediate partition
$\mu$ such that $\lambda<\mu<\sh(S)$ where $\mu=(\lambda_1,\ldots,\lambda_{i+k-2},
\lambda_{i+k-1}-1,\lambda_{i+k}+1,\ldots).$
\par
However, there is no tableau $P$ such that $\sh(P)=\mu$ and $T\gos P\gos S.$
Indeed, $\pi_{1,n-1}(T)=\pi_{1,n-1}(S).$
Assume that there exists $ P$ corresponding to $D_\mu$ such that
$T\gos P \gos S.$ By proposition 4.6.4  we have
$D^{(1,n-1)}_T \leq D^{(1,n-1)}_P \leq D^{(1,n-1)}_S.$
Now $D^{(1,n-1)}_T=D^{(1,n-1)}_S$ so $D^{(1,n-1)}_P=D^{(1,n-1)}_T$
and $D^{(1,n-1)}_P=D^{(1,n-1)}_S.$
\bk
By the first equality $n\in <P^i>$ and by the second equality  $n\in <P^j>$
Since $j >i$ this gives a contradiction.
\QED
In part 2(b) of the proof  we obtain a ``jump'' of length
2 that is $S$ in the cover of $T$ such that $\sh(S)$ is not in the cover $\sh(T)$
however, there exists the unique $\mu$ such that $\mu$ is in the cover of $\sh(T)$
and $\sh(S)$ is in the cover $\mu.$

For example
 $$T = \vcenter{
\halign{& \hfill#\hfill
\tabskip4pt\cr
\multispan{9}{\hrulefill}\cr
\ssa
\vb & 1 &  &  2 & &  4 & \vb &  9 & \ts\vb \cr
\vsa
&&&&&&\multispan{3}{\hrulefill}\cr
\ssa
\vb & 3 &  & 5  & & 8 & \ts\vb \cr
\vsa
&&\multispan{5}{\hrulefill}\cr
\ssa
\vb & 6 & \ts\vb \cr
\vsa
\multispan{0}{\hrulefill}\cr
\ssa
\vb & 7 & \ts\vb \cr
\vsa
\multispan{3}{\hrulefill}\cr}} \ \ \ \
 S =\vcenter{
\halign{& \hfill#\hfill
\tabskip4pt\cr
\multispan{7}{\hrulefill}\cr
\ssa
\vb & 1 &  &  2 & &  4  & \ts\vb \cr
\vsa
&&&&&&\multispan{0}{\hrulefill}\cr
\ssa
\vb & 3 &  & 5  & & 8 & \ts\vb \cr
\vsa
&&\multispan{5}{\hrulefill}\cr
\ssa
\vb & 6 & \vb & 9 & \ts\vb \cr
\vsa
&&\multispan{3}{\hrulefill}\cr
\ssa
\vb & 7 & \ts\vb \cr
\vsa
\multispan{3}{\hrulefill}\cr}} $$
Note that in the example \ref{6.7} we also have a ``jump'' of length 2.

The interesting question is what is the maximal possible length of a
``jump'', i.e the maximal possible length of the chain between $\sh(T)$
and $\sh(S)$  where $S$ is in the cover of $T.$

\subsection{}
\label{7.4} The same ``non-smoothness'' seems to be in charge of the fact that neither
projection, nor injection preserves the cover.
Let us provide the corresponding examples.

We begin with the projection. From our previous discussion it is obvious that
we always have $S$ in the cover of $T$ such that $\pi_{\sr 1, n-\sr 1}(S)=\pi_{\sr 1, n-\sr 1}(T).$
Now we show that there are cases when  $S$ is in the cover of $T$ and there exists $P$
such that $\pi_{\sr 1, n-\sr 1}(T)\gos P\gos \pi_{\sr 1, n-\sr 1}(S).$ The first such example
occurs in $\gs\gl_{\sr 5}.$ Consider
$$T=\vcenter{
\halign{& \hfill#\hfill
\tabskip4pt\cr
\multispan{7}{\hrulefill}\cr
\ssa
\vb  & 1 & & 2 & & 4 &\ts\vb\cr
\vsa
&&&&\multispan{3}{\hrulefill}\cr
\ssa
\vb  & 3 & & 5 &\ts\vb\cr
\vsa
\multispan{5}{\hrulefill}\cr}}\qquad
   \Cscr_T = \left\{
   \begin{array}{c}
s_{\sr 4}s_{\sr 2},\ s_{\sr 4}s_{\sr 1}s_{\sr 2}, \ s_{\sr 3}s_{\sr 4}s_{\sr 2},\\
s_{\sr 3}s_{\sr 4}s_{\sr 1}s_{\sr 2}, \ s_{\sr 2}s_{\sr 3}s_{\sr 4}s_{\sr 1}s_{\sr 2} \\
\end{array} \right\}$$
and
$$
S =
\vcenter{
\halign{& \hfill#\hfill
\tabskip4pt\cr
\multispan{5}{\hrulefill}\cr
\ssa
\vb  & 1 & & 4 &\ts\vb\cr
\vsa
\multispan{0}{\hrulefill}\cr
\ssa
\vb  & 2 & & 5 &\ts\vb\cr
\vsa
&&\multispan{3}{\hrulefill}\cr
\ssa
\vb  & 3 &\ts\vb\cr
\vsa
\multispan{3}{\hrulefill}\cr}}\qquad
   \Cscr_S=\left\{
   \begin{array}{c}
s_{\sr 4}s_{\sr 2}s_{\sr 1}s_{\sr 2}, \ s_{\sr 3}s_{\sr 4}s_{\sr 2}s_{\sr 1}s_{\sr 2}, \
s_{\sr 2}s_{\sr 3}s_{\sr 4}s_{\sr 2}s_{\sr 1}s_{\sr 2},\\
s_{\sr 3}s_{\sr 4}s_{\sr 3}s_{\sr 2}s_{\sr 1}s_{\sr 2}, \
s_{\sr 2}s_{\sr 3}s_{\sr 4}s_{\sr 3}s_{\sr 2}s_{\sr 1}s_{\sr 2} \\
\end{array}\right\}$$
One has $s_{\sr 4}s_{\sr 1}s_{\sr 2}\dos s_{\sr 2}s_{\sr 4}s_{\sr 1}s_{\sr 2},$ hence $T\gos S$.
As well $\pi_{\sr 2,5}(T)=\pi_{\sr 2,5}(S)$ so that $S$ is in the cover of $T.$ However,
$$\pi_{\sr 1,4}(T)=
\vcenter{
\halign{& \hfill#\hfill
\tabskip4pt\cr
\multispan{7}{\hrulefill}\cr
\ssa
\vb & 1 &  & 2 &  & 4 &\ts\vb\cr
\vsa
&&\multispan{5}{\hrulefill}\cr
\ssa
\vb & 3  &\ts\vb\cr
\vsa
\multispan{3}{\hrulefill}\cr}},
\qquad
\pi_{\sr 1,4}(S)=
\vcenter{
\halign{& \hfill#\hfill
\tabskip4pt\cr
\multispan{5}{\hrulefill}\cr
\ssa
\vb & 1 & & 4 &\ts\vb\cr
\vsa
&&\multispan{3}{\hrulefill}\cr
\ssa
\vb & 2  &\ts\vb\cr
\vsa
\multispan{0}{\hrulefill}\cr
\ssa
\vb  & 3 &\ts\vb\cr
\vsa
\multispan{3}{\hrulefill}\cr}}.$$
Note that $\pi_{\sr 1,4}(T)\dos P\dos \pi_{\sr 1,4}(S)$   where
$$P=\vcenter{
\halign{& \hfill#\hfill \tabskip4pt\cr \multispan{5}{\hrulefill}\cr
\ssa \vb & 1 &  & 2 & \ts\vb\cr \vsa &&&&\cr \ssa \vb & 3
&&4&\ts\vb\cr \vsa \multispan{5}{\hrulefill}\cr}}.$$ Thus, $\pi_{\sr
1,4}(S)$ is not in the cover of $\pi_{\sr 1,4}(T).$
\subsection{}
\label{7.5}
As for injections, they do not preserve the cover even in the most trivial case.
One has that $S$ is in the cover of $T$ where
$$T=\vcenter{
\halign{& \hfill#\hfill
\tabskip4pt\cr
\multispan{5}{\hrulefill}\cr
\ssa
\vb & 1 &  & 2 & \ts\vb\cr
\vsa
&&\multispan{3}{\hrulefill}\cr
\ssa
\vb & 3  &\ts\vb\cr
\vsa
\multispan{3}{\hrulefill}\cr}}\qquad {\rm and}\qquad
S=\vcenter{
\halign{& \hfill#\hfill
\tabskip4pt\cr
\multispan{3}{\hrulefill}\cr
\ssa
\vb & 1&\ts\vb\cr
\vsa
&&\cr
\ssa
\vb & 2  &\ts\vb\cr
\vsa
&&\cr
\ssa
\vb & 3  &\ts\vb\cr
\vsa
\multispan{3}{\hrulefill}\cr}}$$
However,
$$
(T\Da 4)=
\vcenter{
\halign{& \hfill#\hfill
\tabskip4pt\cr
\multispan{7}{\hrulefill}\cr
\ssa
\vb & 1 &  & 2 &  & 4 &\ts\vb\cr
\vsa
&&\multispan{5}{\hrulefill}\cr
\ssa
\vb & 3  &\ts\vb\cr
\vsa
\multispan{3}{\hrulefill}\cr}}
\qquad{\rm and}\qquad
(S\Da 4)=
\vcenter{
\halign{& \hfill#\hfill
\tabskip4pt\cr
\multispan{5}{\hrulefill}\cr
\ssa
\vb & 1 & & 4 &\ts\vb\cr
\vsa
&&\multispan{3}{\hrulefill}\cr
\ssa
\vb & 2  &\ts\vb\cr
\vsa
\multispan{0}{\hrulefill}\cr
\ssa
\vb  & 3 &\ts\vb\cr
\vsa
\multispan{3}{\hrulefill}\cr}},$$
so that $(S\Da 4)$ is not in the cover of $(T\Da 4)$ exactly as in \ref{7.4}.

\subsection{}\label{7.6}
Let us finish with a very simple lemma showing that $\Tscr_{\al,\be}$ preserves
the cover. We give its one line proof for the completeness.
\begin{lemma}
Let $T,S\in D_{\al,\be}.$ Then $S$ is in the cover of $T$ iff $\Tscr_{\al,\be}(S)$
is in the cover of $\Tscr_{\al,\be}(T).$
\end{lemma}
\Pf
Indeed, by the symmetry it is enough to show that if $S$ is in the cover of $T$
then $\Tscr_{\al,\be}(S)$ is in the cover of $\Tscr_{\al,\be}(T).$

Assume that $S$ is in the cover of $T$, but $\Tscr_{\al,\be}(S)$ is not in the cover of
$\Tscr_{\al,\be}(T).$ Then there exists $P$ such that $\Tscr_{\al,\be}(T)\gos P\gos \Tscr_{\al,\be}(S).$
The fact that $\Tscr_{\al,\be}(T),\Tscr_{\al,\be}(S)\in D_{\be,\al}$ forces $P\in D_{\be,\al}$
so that by applying $\Tscr_{\be,\al}$ to all three tableaux we get $T\gos \Tscr_{\be,\al}(P)\gos S$
which contradicts $S$ being in the cover of $T.$
\QED
The facts that for a  Richardson tableau $T$  one has $S\dg T$ iff $S\chg T$ and for both $\dos$ and $\cho$
one has $T>S$ iff $T^\dagger< S^\dagger$ leaded us to conjecture that both the geometric and
the algebraic orders can be obtained only by considering the order on Richardson components together
with operations $\Tscr_{\al,\be}$ and transposition. However, these 3 operations are enough to construct
geometric order only for $n\leq 7.$ The computations show that for $n=8$ there is a pair $S,T$
such that $S$ is in the cover of $T$ in Duflo order, however, $S$ cannot be obtained as an element
of the cover of $T$ with the help of our 3 operations. In this connection the fact that for $n=9$
these 3 operations again give us the full picture seems to be even more peculiar. Of course, for $n\geq 10$
where the geometric order does not coincide with Duflo-Vogan order anymore, these 3 operations cannot give
the full picture.
\bigskip
\parno
\centerline{ INDEX OF NOTATION}

\parno
\begin{tabular}{llll}
\ref{1.2}& $\bG,\ \gog,\ U(\gog),\ \Oscr,\ \gn,\ \gn^-,\ \gh,\
\Vscr_w;$&
         \ref{2.4}& $\bP_{\al},\ \Iscr,\ \bP_{\Iscr},\
       \bM_{\Iscr},\ \bL_{\Iscr},\ \gp_{\sr\Iscr},\ \gm_{\sr \Iscr},\
        \gl_{\sr\Iscr},\bB_{\Iscr},$\\
\ref{1.5}& $R,\ R^+,\ \Pi,\ \rho,\ L_w;$&
      &$\gn_{\sr \Iscr},\ \pi_{\sr\Iscr},\ W_\Iscr,\ F_\Iscr,\ R^+_\Iscr,\ w_{\sr \Iscr},\
        \Cscr^\Iscr_{w_\Iscr},\ \Vscr_{w_\Iscr}^\Iscr;$\\
\ref{1.6}& $\bT_n,\ I_T,\ \Vscr_T,\ \ao,\ \go;$&
  \ref{2.6}&$(T)_{i,j},\, [w,y],\, \ov w_a,\,  \ov T_a,\, \Iscr_{i,j},\, \pi_{i,j},\, \bS_{i,j},\
            T^{\dagger},\, w_o;$\\
\ref{1.7}& $\Oscr_{\Vscr};$&
    \ref{3.1}& $ D_\lambda\leq D_\mu;$\\
\ref{1.8}& $\prec,\ \preceq;$&
   \ref{3.4}& $\gn_k,\ \theta(u),\ \phi(T),\ \nu_{\sr T},\ \Vscr_T;$\\
\ref{1.9}& $I_w,\ \ell(w), \dor,\ \cho;$&
   \ref{3.9}& $\varphi(T);$\\
\ref{1.10}& $\Pi,\ \tau(T),\ D_{\al,\be},\ \Tscr_{\al,\be};$&
    \ref{4.1}& $\Xscr_0,\ \gb,\ M_w,\ L_w,\ I_w;$\\
\ref{1.11}& $\dvo,\ \chvo;$&
    \ref{4.2}& $\Cscr_L,\ I_{\Cscr},\ \Cscr_D,\ \ao;$\\
\ref{2.1}& $\bB,\, X_{\al},\, w(X_\al),\,
       \ov{\bB(\gn\cap^w\gn)},\, \Oscr_w,\, \Vscr_w;$&
       \ref{4.3}& $S(w),\ \tau(w),\ \tau(I),\ \tau(\Vscr);$\\
\ref{2.2}& $\Cscr_w,\, D_{\lambda},\, \bD_n,\, \sh(T),\, \bT_n,\,
      \bT_{\lambda},\, $&
       \ref{4.5}& $\rho_{\sr\Iscr},\ M_w^\Iscr,\ L_w^\Iscr,\ I_w^\Iscr$\\
 & $P(w),\, Q(w),\, J(u),\, \Oscr_u,\, \Nscr,\, \Oscr_{\lambda},\ D_u$&
       \ref{5.1}& $D_{\al,\beta},\ \Tscr_{\al,\be}$\\
\ref{2.3}& $T\dor S$&
       \ref{6.7}& $\dvo,\ \chvo$\\

\end{tabular}

\end{document}